\documentclass[runningheads]{llncs}
\usepackage[T1]{fontenc}
\usepackage{graphicx}
\usepackage{booktabs}
\usepackage[misc]{ifsym}

\usepackage{mwe}

\usepackage[misc]{ifsym}
\usepackage{amssymb}
\usepackage{amsmath}
\usepackage{cases}

\usepackage{xr-hyper}
\usepackage{cleveref}

\usepackage[square,sort,comma,numbers]{natbib}
\usepackage{tikz}

\usepackage{float}

\newcommand{\setJK}{{\cal{J}}_K}
\newcommand{\setJ}{{\cal{J}}_k}
\newcommand{\setR}{\mathbb{R}}
\newcommand{\setX}{\mathcal{D}}
\newcommand{\setJKbar}{{\bar{\cal{J}}}_K}

\usepackage{cases}
\usepackage{wrapfig}
\usepackage[table]{xcolor}
\definecolor{lblue}{rgb}{0.94, 0.94, 0.94}

\usepackage{multirow}
\usepackage{subcaption}

\usepackage[algo2e,boxed,lined]{algorithm2e}
\toctitle{Fast SDP certification of neural networks : towards large multi-class datasets}

\author{\Letter{Margot Boyer}\inst{1}, Clément Rambour\inst{2}, Zacharie Alès\inst{1,3}, Amélie Lambert\inst{1}}

\tocauthor{\Letter{M. Boyer}, C. Rambour, Z. Alès, A. Lambert}

\begin{document}

\title{Fast SDP certification of neural networks : towards large multi-class datasets}

\titlerunning{Fast SDP certification of neural networks : towards large multi-class datasets}

\authorrunning{M. Boyer et al.}



\institute{Cedric, Conservatoire National des Arts et Métiers, F-75003 Paris, France\email{\{margot.boyer, amelie.lambert\}@lecnam.net}
\and
Sorbonne Université, CNRS, Institut des Systèmes Intelligents et de Robotique, ISIR, F-75005 Paris, France \email{clement.rambour@isir.upmc.fr}
\and
Institut Polytechnique de Paris, ENSTA, Unité de Mathématiques Appliquées, F-91120 Palaiseau, France \email{zacharie.ales@ensta.fr}
}

\maketitle              

\begin{abstract}
 We present a new quadratic model for the certification problem in adversarial robustness, which simultaneously accounts for all possible target classes. Building on this model, we propose a novel semidefinite programming (SDP) relaxation for incomplete verification. A key advantage of our approach is that it certifies robustness in a single optimization, avoiding the need for a separate resolution per class. This yields a significant computational speed-up and enables scalability to large datasets with many classes. To further improve efficiency, we also propose an effective pruning strategy of active neurons, thus reducing the problem dimensionality and accelerating convergence.

\keywords{Adversarial robustness \and Certification \and SDP optimization.}
\end{abstract}

\section{Introduction}
Deep Neural Networks (DNNs) have achieved remarkable success and are widely implemented in various domains, including computer vision and natural language processing. The rapid adoption of DNNs has often prioritized efficiency and automation, sometimes at the expense of safety considerations. 

The research community has extensively studied various aspects of robustness, including out-of-distribution generalization, robustness to data corruption, and resistance to adversarial attacks. In particular, DNN have been proven vulnerable to adversarial attacks~\citep{goodfellow_explaining_2015}, where malicious actors exploit the complexity of these models to generate examples that deceive classifiers. This issue has raised concerns in many critical domains of application of neural networks, like autonomous vehicles or robotics, where adversarial attacks could be a means for malicious purposes.

An \emph{adversarial attack} consists in solving a constrained optimization problem to find an adversarial example for a given input $x$, \textit{i.e.}, a  data in the neighborhood of $x$ which is classified differently by the DNN. These attacks represent a significant threat, particularly when the attacker has knowledge of the model architecture and parameters. In response, two main approaches have emerged to enhance the robustness of DNNs against such attacks: adversarial training and certified defenses.
Adversarial training methods aim to improve robustness by performing adversarial augmentations. 
While these methods do offer increased resilience, they are not foolproof and can still be vulnerable to sophisticated attacks. 

On the other hand, certified defenses provide mathematical guarantees of robustness against adversarial attacks. 

 The certification problem for neural networks with ReLU activation functions is NP-complete~\citep{katz_reluplex_2017}. This inherent complexity implies that providing a \emph{complete certification} requires substantial computational effort and remains limited in scalability. Many approaches solve combinatorial models to assess the stability of DNN predictions around each data. Several Mixed-Integer Programming (MIP) formulations were introduced to provide formal proofs for small ReLU DNNs ~\citep{tjeng_evaluating_2019,fischetti_deep_2018, dsouza_maximum_2017} but remain intractable for medium to large-scale problems.

Computing a non-negative lower bound is sufficient to certify that no adversarial example exists for a given target class. Thus, in order to speed up the certification, many approaches solve a relaxation of the original certification problem.
 In this paper, we focus on \emph{incomplete} verifiers that provide lower bounds on the certification problem: a positive bound guarantees robustness while a negative bound is inconclusive. They constitute a compromise between efficiency and scalability, aiming to achieve the highest possible lower bound within a time limit. 
 
 Most of the incomplete certification methods are based on quadratic optimization formulations in which the ReLU is expressed as a quadratic non-convex equality from which a linear relaxation is computed~\citep{wong_provable_2018}. 
Despite bringing promising results for certification, current approaches using SDP relaxations offer limited scalability in particular when certifying mid-to-large scale datasets composed of multiple classes. 
Indeed, these approaches  are targeted ones, i.e. each combinatorial model tests if there exists an adversarial attack for one data and one \textit{target} class. Thus, formally certifying a single data point requires looping over all possible target classes. This requirement quickly becomes cumbersome as modern datasets such as ImageNet-1k or ImageNet-21k propose hundreds or thousands of classes. Furthermore, as each neuron brings its own set of constraints in the optimization problem, current SDP approaches struggle for deep networks. 

\paragraph{\textbf{Contributions}} Our contributions are threefold: 
\begin{itemize}
    \item   To deal with the aforementioned limitations, we introduce a new model for the certification problem that is based on an untargeted quadratic formulation $(QP_U)$. This new approach allows us to certify each data by solving a \emph{single} optimization problem and significantly speeds up the certification process. A key advantage of our formulation is that it preserves the non-negativity condition, certifying the data whenever a non-negative lower bound of $(QP_U)$ is obtained. \\[1ex]
    \item We further introduce valid quadratic inequalities that tighten the bound of the relaxed problem. Finally, to scale up the certification, we propose an efficient pruning strategy able to remove \emph{all stable} neurons from the optimization constraints. This allows us to reduce the solution time. An interesting result is that this pruning strategy is generic and can be applied to other SDP relaxations. \\[1ex]
    \item Finally, we present computational results demonstrating the efficiency of our methods against state-of-the-art approaches. 
\end{itemize}

\section{Related Work}
\paragraph{\textbf{Certification problem}} 
Complete certification methods aim to provide definitive guarantees about the absence of adversarial examples within a given input region. Seminal works include approaches based on Mixed-Integer Linear Programming (MILP), which model ReLU activations through integer constraints, enabling exact reasoning over the network’s activations. The MILP formulations proposed by \citep{fischetti_deep_2018,tjeng_evaluating_2019} demonstrate formal verification for small-sized ReLU networks. Similarly, satisfiability modulo theories solvers have been employed to provide sound and complete verification \citep{ehlers_formal_2017,katz_reluplex_2017}. 
However, a key limitation of complete verifiers remains their limited scalability to deep networks or high-dimensional datasets.
\paragraph{\textbf{Incomplete verifiers}} 
aim to compute tight lower bounds on neural network robustness, providing formal certification guarantees whenever these bounds are sufficiently strong.
A positive bound confirms robustness, while a negative one remains inconclusive, making these methods an interesting compromise between scalability and theoretical soundness. Incomplete verifiers are divided into a wide variety of approaches including convex relaxations via duality~\citep{wong_provable_2018,Gowal_ijcai2019}, linear bounding of ReLU  activations~\citep{weng_towards_2018}, or discretized input space exploration~\citep{huang2017safety}. However, most incomplete verifiers may yield conservative bounds, and even well-optimized linear relaxations can fail to produce tight lower bounds over the objective.
\paragraph{\textbf{Semidefinite Programming (SDP) relaxations}}
have emerged as a promising class of incomplete verifiers for neural network verification as they produce sharper lower bounds than traditional linear programming approaches. While SDP methods are computationally more intensive, foundational works~\citep{raghunathan_semidenite_2018, zhang_tightness_2020,dathathri_enabling_2018,chiu_tight_2023,lan_iteratively_2023} have demonstrated empirical tightness compared to LP relaxations. This was further advanced by integrating geometric constraints such as triangle relaxations~\citep{batten_efficient_2021} and Reformulation-Linearization Technique (RLT) cuts~\citep{lan_tight_2022}, which refine the feasible region for ReLU -activated networks. However, the relaxation becomes increasingly loose as the network depth grows, and solving SDPs for deep networks often results in scalability challenges. This phenomenon is exacerbated when dealing with multiple classes as one SDP relaxation needs to be computed \emph{for each target class} to achieve certification.

\section{Preliminaries}
Deep Neural Networks (DNNs) considered in this work are non-linear functions that map the input set to a measurable label set. They are described as successive layers given by the composition of a linear and a non-linear transformation. 
  
Each layer $k$ contains $n_k$ neurons, indexed by $\setJ =\{1,...,n_k\}$. The output $z_{k+1}\in\setR^{n_{k+1}}$ of every layer $k \in \left[K-2\right]$ (\textit{i.e.} $\{0,...,K-2\}$) is computed by a ReLU activation function:  $z_{k+1} = ReLU(W_{k+1} z_k + b_{k+1})$, where $b_k \in \mathbb{R}^{n_k}$, $W_k \in \mathbb{R}^{n_{k} \times n_{k-1}}$ are the learned parameters of the network. Given a finite labeled dataset $\mathcal D = \{x_i,y_i\} $, 
the predicted class is given by $y^* = \textrm{argmax}_{j \in \setJK} z^j_K$ for data $x=z_0$ where $z^j_K$ is the $j$-th component of vector $z_K=W_{K} z_{K-1} + b_{K}$. 

For a given $\epsilon > 0$, the certification task verifies that for each data $( {x},{y})$ and all $z_0\in \mathcal{B}_{\epsilon}({x})$ (the $\infty$-norm ball centered at ${x}$ and radius $\epsilon$) the DNN correctly predicts the class ${y}$.

Defining $\setJKbar = \setJK\backslash\{ y\}$, the set of all possible targets for a given sample, $W^j_K$, the $j^{th}$ row of matrix $W_K$, and $\mathcal D^+$, the set of well-classified data and their labels, we formally define robustness as follows: 

\begin{proposition}[Targeted Robustness]
    For a data point $(x,y)\in\mathcal D^+$,  a target class $j\in\setJKbar$ and $\epsilon >0$,  a neural network is $(\epsilon,j)-$robust in $x$, if :$$\min\limits_{z_0 \in \mathcal{B}_{\epsilon}(x)} z_K^{y} - z_K^{j} \geq 0 $$
     \label{propTarget}
\end{proposition}

\begin{proposition}[Full robustness]
    For $\epsilon >0$, a neural network is $\epsilon-$robust if  for all  $(x,y) \in  \mathcal D^+$, \begin{equation}\min\limits_{j \in \setJKbar} \,\, \min\limits_{z_0 \in \mathcal{B}_{\epsilon}(x)} z_K^{y} - z_K^{j} \geq 0 \label{prop:2} \end{equation}
    \label{propFull}
\end{proposition}

Our objective in this paper is to determine whether the DNN satisfies the robustness condition defined in Proposition~\ref{propFull}. More formally, considering a data $(x,y) \in \mathcal D^+$, we consider the following optimization problem $(\textrm{Cert})$ defined for all $j \in \setJK$ :

\begin{align}
    (\mathrm{Cert}) \quad 
    \min_{z_0} \;& 
    \big(W^y_K z_{K-1} + b^y_K\big)
    -
    \big(W^j_K z_{K-1} + b^j_K\big) \\[0.5ex]
    \text{s.t.} & z_{k+1}=ReLU(W_{k+1}z_k+b_{k+1}), \; \forall k\in[K-2] \label{ReLU} \\[0.3ex] 
    & x - \epsilon \le z_0 \le x + \epsilon \label{ball}
\end{align}
where Constraints~(\ref{ReLU}) fix the output of the hidden layers and Constraint~(\ref{ball}) ensures that $x$ belongs to $\mathcal{B}_{\epsilon}(x)$. The objective is the difference between the logit of the true class $y$ and the target class $j$. 

Solving $(\textrm{Cert})$ to global optimality is hard due to the non-convexity of Constraints~(\ref{ReLU}). However, by denoting $v(\textrm{Cert})$ the optimal value of $(\textrm{Cert})$, Proposition~\ref{propFull} is reached when $v(\textrm{Cert}) \geq 0$, for all $j \in \setJK$ and $(x,y) \in \setX^+$. Thus, it is sufficient to compute a non-negative lower bound of $(\textrm{Cert})$ for all $j \in \setJK$ to ensure full robustness. 

  A quadratic formulation of $(\textrm{Cert})$ was introduced in~\citep{raghunathan_semidenite_2018}, obtaining the following targeted formulation:
\begin{numcases}{(QP_{T}^j)}
\underset{z_0}{\min} \; \big(W^y_{K} z_{K-1} + b^y_{K}\big) - \big( W^j_{K} z_{K-1} + b^j_{K}\big)  \nonumber \\
         \mbox{s.t.}
         \quad z_{k+1} \geq 0, \quad z_{k+1} \geq W_{k+1} z_k + b_{k+1}   \label{ReLURef1}, \; \;  \forall k \in [K-2], \\
         \quad z_{k+1} \odot (z_{k+1} -  W_{k+1} z_k - b_{k+1}) = 0, \; \;  \forall k \in [K-2], \label{ReLURef3} \\
         \quad z_{k} \odot z_{k}- (L_{k} + U_{k}) \odot z_k + U_{k} \odot L_{k} \leq 0, \; \; \forall k \in [K-1] \label{boundRef}
    \end{numcases}
   where $\odot$ denotes the element-wise product, $L_k$ and $U_k$ are lower and upper bounds over the preactivation vector of layer $k$.  Constraints~(\ref{ReLURef1}) combined with Constraints~(\ref{ReLURef3}) are equivalent to Constraints~(\ref{ReLU}). 
    Constraints~(\ref{boundRef}) can be rewritten as $(U_k - z_k) \odot (z_k - L_k) \geq 0$, which enforces $L_k \leq z_k \leq U_k$ when $L_k \leq U_k$. For $k=0$, this is a quadratic equivalent to Constraints~(\ref{ball}). Note that there exist efficient methods to propagate bounds across the network starting from the bounds of the input layer (e.g., $L_0 = x - \epsilon$ and  $U_0 = x + \epsilon$ with the $\infty$ norm)~\citep{wang_beta-crown_2021} giving bounds $L_k$ and $U_k$ on the preactivation vector for all layers $k$.

   \begin{proposition}[Target-positivity property]
     If the optimal value of the targeted quadratic formulation $v(QP_{T}^j)$ is non-negative, the DNN satisfies Proposition~\ref{propTarget}.
     \label{posTarget}
   \end{proposition}

      Due to the non-convexity of Constraints~(\ref{ReLURef3}) and~(\ref{boundRef}), solving formulation $(QP_{T}^j)$ to global optimality is impractical even for small-sized DNNs. However,  Proposition~\ref{posTarget} ensures that the development of suitable relaxations can be sufficient to certify the robustness. In particular, using semi-definite relaxations for quadratic programming was widely studied~\citep{Ans09}. Let $P = \begin{bmatrix} 1 & z \end{bmatrix} \begin{bmatrix} 1 & z \end{bmatrix}^T$ be the matrix that collects all the linear and quadratic terms in $(QP_{T}^j)$. Then, the semi-definite relaxation of targeted problem $(QP_{T}^j)$ has the form:
\begin{numcases}{(SDP_{T}^j-IP)}
   \underset{P}{\min} \;\; \big(W^y_{K} P[z_{K-1}] + b^y_{K}\big) - \big(W^j_{K} P[z_{K-1}] + b^j_{K}\big)      \nonumber \\
   \mbox{s.t.}
   \quad P[z_{k+1}] \geq 0, \quad P[z_{k+1}] \geq W_{k+1} P[z_k] + b_{k+1}, \; \;  \label{ReLURef1_SDP} \\
     \quad \textrm{diag} \big(P[z_{k+1}z_{k+1}^{\top}] -  W_{k+1} P[z_k z_{k+1}^{\top}] \big) = b_{k+1} \odot P[z_{k+1}], \; \;   \label{ReLURef3_SDP} \\
     \quad\quad\quad\quad\quad\quad\quad\quad\quad\quad\quad\quad\quad\quad\quad\quad\; \forall k \in [K-2] \\
     \quad \textrm{diag}\big(P[z_{k}z_k^{\top}]\big) - (L_{k} + U_{k}) \odot P[z_k] + U_{k} \odot L_{k} \leq 0 \label{boundRef_SDP} \\\quad\quad\quad\quad\quad\quad\quad\quad\quad\quad\quad\quad\quad\quad\quad\quad\; \forall k \in [K-1] \\
    \quad P \succeq 0, P[1] = 1  \label{SDP} 
 \end{numcases} 
\noindent From its definition, each element of $P$ is related to a given term in (Cert) and the symbolic indexing $P[.]$ is used to index the vectors of elements of matrix $P$. Constraints~(\ref{ReLURef1_SDP})--(\ref{boundRef_SDP}) correspond to the linearization of Constraints~(\ref{ReLURef1})--(\ref{boundRef}) by the matrix variable $P$.

 
 Note that matrix $P$ exhibits a block diagonal structure.
 By leveraging chordal decomposition techniques~\citep{vandenberghe_chordal_2015}, as outlined in~\citep{batten_efficient_2021}, $P$ can be decomposed into multiple submatrix variables. Specifically, the decomposition yields $K-1$ matrix variables, each associated to two consecutive layers $k$ and $k+1$ (for $k \in [0,K-2]$): $P_k = \begin{bmatrix}1 & z_k & z_{k+1}\end{bmatrix} \begin{bmatrix}1 & z_k & z_{k+1} \end{bmatrix}^{\top}$. This decomposition allows to deal with multiple modest-sized SDP matrices rather than $P$. This change is expressed in the previous constraints (\ref{ReLURef1_SDP})--(\ref{boundRef_SDP}) by injecting these matrices as: 

\begin{numcases}{}
   \quad P_k[z_{k+1}] \geq 0, \quad P_k[z_{k+1}] \geq W_{k+1} P_k[z_k] + b_{k+1}, \;\;  \forall k \in [K-2] \, \label{ReLURef1_SDP_k} \\
     \quad \textrm{diag} \big(P_k[z_{k+1}z_{k+1}^{\top}] -  W_{k+1} P_k[z_k z_{k+1}^{\top}] \big) = b_{k+1} \odot P_k[z_{k+1}], \;\;  \forall k \in [K-2] \,\label{ReLURef3_SDP_k} \\     
     \quad \textrm{diag}\big(P_k[z_{k}z_k^{\top}]\big) - (L_{k} + U_{k}) \odot P_k[z_k] + U_{k} \odot L_{k} \leq 0, \;\; \forall k \in [K-1]  \label{boundRef_SDP_k} \\
    \quad P_k[z_{k+1}] \leq A_{k+1} P_k[z_k] + B_{k+1}, \;\;  \forall k \in [K-2] \, \label{tri_SDP} 
\end{numcases}

The triangular constraint~(\ref{tri_SDP}) introduced in~\citep{ehlers_formal_2017} tightens the upper bounds of a neuron $j$ of layer $k$ according to its activation status with $A_{k} = l_{k} \odot W_{k}$, $B_{k} = l_{k} \odot (b_{k} - L_{k}) + \text{ReLU}(L_{k})$.

The variable $l_{k} = \frac{\text{ReLU}(U_k) - \text{ReLU}(L_k)}{U_k - L_k}$ indicates whether a neuron is stable active ($l_{k} = 1$), stable inactive ($l_k =0$), or unstable ($0 < l_{k} < 1$). In the case of a \emph{stable active} neuron, \textit{i.e}. $L_k^{j} \geq 0$ then $z^j_{k} \leq W_k^{j} z_{k-1} + b_k^{j}$. If it is \emph{stable inactive}, \textit{i.e.} $U_k^{j} \leq 0$ then $z^j_{k} \leq 0$. In the \emph{unstable} case, constraint~(\ref{tri_SDP}) reduces to $z^j_{k} \leq \frac{U_k^j}{U_k^j - L_k^j} (W_k^{j} z_{k-1} + b_k^{j})+ \frac{U_k^j}{U_k^j - L_k^j}(b_k^j - L_k^j)$.  
 
 It has been further improved in~\citep{lan_tight_2022} with the addition of  RLT (Reformulation Linearization Technique)  cuts, giving the following $(SDP^j_T)$ problem:
 \begin{numcases}{(SDP^j_T)}
   \min \big(W^y_{K} P_{K-2}[z_{K-1}] + b^y_{K}\big) - \big( W^j_{K} P_{K-2}[z_{K-1}] + b^j_{K}\big)     \nonumber \\
   \mbox{s.t. }  (\ref{ReLURef1_SDP_k})-(\ref{tri_SDP}) \nonumber\\ 
    \quad P_{k}[(1 \, z_{k+1}) (1 \, z_{k+1})^{\top}] = P_{k+1}[(1 \, z_{k+1}) (1 \, z_{k+1})^{\top}]  \label{stab_SDP} \\ 
     \quad  RLT(p)  \label{RLT_p_SDP} \\ 
    \quad P_k \succeq 0, P_k[1] = 1, k \in [K-2]\label{SDP_rec} 
\end{numcases} 
\noindent  
where the minimization over $P$ has been omitted for clarity of notation. Constraints~(\ref{stab_SDP}) ensure the coherence of the variables across consecutive matrices $P_k$ and $P_{k+1}$.  For $L_k \leq z_k \leq U_k$ $\forall k \in [K-1]$, Constraints~(\ref{RLT_p_SDP}) are RLT cuts~\citep{SheAda90}. They result from the product of linear valid inequalities (triangular constraint, bound constraints) to obtain quadratic inequalities. A percentage $p$ of these RLT cuts is carefully chosen by a heuristic (see Appendix~\cref{supp_seq:app_rlt} ).

Using $(SDP_{T}^j)$ to certify a DNN requires solving one SDP for each data $(x,y) \in \setX^+$ and each possible target ($j \in \setJKbar$). This leads to two significant drawbacks. First, the certification process fails to scale with the number of classes. Second, since solving multiple SDP relaxations for \textbf{each} data point is computationally demanding, the number of cuts must be restricted, which in turn weakens the tightness of the resulting bounds on the objective. To answer these limitations, we now introduce a new \emph{untargeted} model that certifies a sample for all targets by solving only one SDP, thus considerably reducing the certification burden.

\section{Method}
\subsection{A new quadratic model for full certification}

\begin{figure*}[t!]
    \centering
    \includegraphics[width=0.95\linewidth]{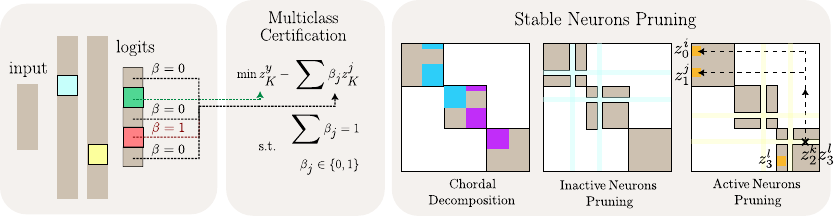}
    \caption{\textbf{Multiclass certification} Our method provides a new quadratic formulation for certifying a DNN across all labels simultaneously. It relies on binary variables $\{\beta_j\}$ indicating the class associated with the worst adversarial example. 
    To further reduce the number of variables, we also propose a pruning strategy that removes both inactive (blue) and active (yellow) neurons. By use of a chordal decomposition of the SDP matrix (in blue and purple), we remove terms required to express the ReLU activations, which are compensated with the introduction of dedicated constraints. For instance, once neuron $z_2^k$ is removed, the quadratic interaction $z_2^k z_3^l$ is no longer represented in the matrix. To address this, we bound the dependencies of $z_3^l$ with respect to neurons from previous layers \textit{i.e.} $z_0^i$ and $z_1^j$.}
    \label{fig:outline}
    
\end{figure*}

To avoid solving $|\setJKbar|$ SDPs for each data $(x,y)\in\mathcal X^+$, we design a new model that directly checks Proposition~\ref{propFull}. As shown in the left part of~\cref{fig:outline}, we thus introduce binary variables $(\beta_j)_{j \in \setJKbar}$ with $\beta_j$ equals $1$ if and only if the worst adversarial example is of class $j$. Thus, the left-hand side of Equation~(\ref{prop:2}) can be  obtained by minimizing $z_K^{y} - \displaystyle{\sum_{j \in \setJKbar} } \beta_j z_K^{j}$.
Using this objective, we define the following quadratic formulation of the \emph{full-robustness} (see Proposition~\ref{propFull}) of a DNN:
\begin{numcases}{(QP_{U})}
\min W^y_{K} z_{K-1} + b^y_{K} - \displaystyle{\sum_{j \in \setJKbar} } \beta_j \big(W^j_{K} z_{K-1} + b^j_{K}\big)  \nonumber \\
         \mbox{s.t.} \,\, (\ref{ReLURef1})-(\ref{boundRef}) \nonumber\\
         \quad \displaystyle{\sum_{j \in \setJKbar} } \beta_j =1 \label{beta_sum1}\\
         \quad  \beta_j \in \{0,1\} & $j \in \setJKbar$  \label{binary} 
\end{numcases}

  Constraint~(\ref{beta_sum1}) ensures that only one binary variable $\beta_j$ will be non-zero. 
 We now prove that $v(QP_U)$ coincides with the lowest value of $v(QP_T^{j})$ over all target classes.
 \begin{theorem}
Given a data $(x,y) \in \setX^+$, $\epsilon > 0$, we have $v(QP_{U}) =\min\limits_{\bar j \in \setJKbar} v(QP_{T}^{\bar j})$ 
\label{prop:equi}
\end{theorem}
From this theorem, we can deduce that the non-negativity of  $v(QP_U)$ ensures the full robustness. Indeed, an optimal value of $(QP_{U})$ will set $\beta_{\bar j}=1$ for the target class $\bar j$ which minimizes $z_K^{y} -  z_K^{\bar j}$. 

 \begin{proposition}[Full positivity property]  
    If $v(QP_{U})$ is non-negative, the DNN satisfies Proposition~\ref{propFull}.   
     \label{posFull}
   \end{proposition}

$(QP_{U})$ has only $|\setJKbar|$ additional binary variables and $1$ more constraint than $(QP_{T}^j)$. Note that one can prune some of these binary variables based on the lower and upper bounds of their logits: a class $j \in \setJKbar$ can be pruned if there exists another class $\tilde{j} \in \setJKbar$ such that $U_K^{j} \leq L_K^{\tilde{j}}$ since the associated $\beta_j$ will be zero in any optimal solution of $(QP_U)$. Indeed, if there exists an adversarial attack, the most damaging one would not target class $j$, but a class in $\setJKbar \backslash \{j\}$. In the following, we assume that the dominated classes are excluded from the set $\setJKbar$~\ref{supp_seq:algos}. This implies $\bigcap\limits_{j\in \setJKbar}  [L_K^{j}, U_K^{j}] \neq \emptyset$.

 Similarly to $(QP_{T}^j)$, the direct solution of $(QP_{U})$ is impractical even for small-sized DNNs. Thus, we build a tight semi-definite relaxation of $(QP_{U})$ that may certify the DNN using Proposition~\ref{posFull}.

\subsection{A tight SDP relaxation}

In order to handle the additional binary variables $\beta_j$ in our SDP relaxation, for $k\in [K-3]$, we use $P_k$ as defined in Constraint~(\ref{SDP}) and
we introduce an additional matrix $P_{K-2} = \begin{bmatrix}
    1 & z_{K-2 } & z_{K-1}& \beta
\end{bmatrix}
\begin{bmatrix}
    1 & z_{K-2 } & z_{K-1}  & \beta
\end{bmatrix}^{\top}$ to linearize the products $\beta_jz_K$.

We build the following SDP relaxation of $(QP_U)$:
\begin{numcases}{(SDP)}
    \min \; W_{K}^{y} \; P_{K-2}[z_{K-1}] + b_{K}^{y}
    - \sum_{j \in \setJKbar}  \big( W_{K}^{j} \; P_{K-2}[\beta_j z_{K-1}] + b_{K}^{j} \big) \nonumber \\
    \text{s.t. }  (\ref{ReLURef1_SDP_k})-(\ref{SDP_rec})  \nonumber \\
    \quad \sum_{j \in \setJKbar} P_{K-2}[\beta_j] = 1  \label{beta_sum1_SDP} \\
    \quad \mathrm{diag}\big(P_{K-2}[\beta \beta^{\top}]\big) = P_{K-2}[\beta] \label{binary_SDP} 
\end{numcases}

To tighten $(SDP)$, we use several cuts. First, to tighten the linearization of products $\beta$ by $z$, we use McCormick cuts~\citep{Cor76} for $k \in \{K-2, K-1\}$ present in matrix $P_{K-2}$, that are defined as follows:
\begin{equation}{ }
 \beta_{j} z_{k} \geq 0, \quad \, \, \beta_{j} z_{k} \geq U_k \beta_{j} + z_k - \beta_{j}, \quad \, \, 
 \beta_{j} z_{k} \leq U_k \beta_{j},\quad  \, \,
 \beta_{j} z_{k} \leq z_k   \label{mc1}
 \end{equation}

Note that logits $z_K^{j}$ do not appear explicitly in our formulation, instead they are represented through their linear combination in the penultimate layer. Similarly, we tighten the products $\beta_{j_1}z^{j_2}$ of the objective function by use of the McCormick envelopes for all $(j_1,j_2) \in \setJK \times \setJK$:
 \begin{equation}
 \beta_{j_1}z^{j_2}_{K} \leq z^{j_2}_{K} - L_{K-1}^{j_2} (\beta_{j_1}-1), \,\, \quad 
  \beta_{j_1}z^{j_2}_{K} \leq U_{K}^{j_2} \beta_{j_1}, \, \, \nonumber
  \end{equation}
   \begin{equation}
  \beta_{j_1}z^{j_2}_{K}\geq z^{j_2}_{K} - U_{K}^{j_2} (\beta_{j_1}-1), \,\, \quad 
  \beta_{j_1}z^{j_2}_{K} \geq L_{K}^{j_2} \beta_{j_1}   \label{mc4-obj}
 \end{equation}
Then, we build 3 new specific families of valid quadratic cuts for $(QP_{U})$. First,  since $\beta$ is a unit vector, a subset of entries can be fixed, reducing the number of terms that require explicit modeling, and we get:
\begin{equation}\beta_{j_1} \beta_{j_2}  =0 \;\;\,\, \forall (j_1,j_2) \in \setJKbar, j_1 \neq j_2 \label{cut:betaibetaj} \end{equation}
 Finally, for all pairs of two distinct adversarial targets $j_1,j_2$, we introduce new inequalities that leverage the specific structure of the certification problem. These constraints allow coupling the variables $\beta_{j_1}$, $\beta_{j_2}$, $z_K^{j_1}$, $z_K^{j_2}$, $\beta_{j_1} z_K^{j_1}$ and $\beta_{j_2} z_K^{j_2}$ using only two constraints. Intuitively, constraint~(\ref{eq:cut2}) encourages the logit of the adversarial target selected by the model to exceed that of other possible targets. 
 \begin{numcases}{}
    \beta_{j_2} z_{K}^{j_2} \leq (1 - \beta_{j_1}) z_{K}^{j_1} + \beta_{j_2} U_{K}^{j_2} - (1 - \beta_{j_1}) L_{K}^{j_1} +\beta_{j_2} (L_K^{j_1}-z_{K}^{j_1})  &  \label{eq:cut1}\\
    \beta_{j_2} z_{K}^{j_2} \geq z_{K}^{j_1} - (1 - \beta_{j_2}) \; U^{j_1}_K& \label{eq:cut2}
\end{numcases}
\begin{proposition}
    Constraints~(\ref{eq:cut1}) and~(\ref{eq:cut2}) are valid inequalities of ($QP_U$).
\end{proposition}

\subsection{Pruning of stable active neurons}

 Formulating the certification problem with the smallest possible model is crucial to accelerating its resolution. With this in mind, recent works have shown~\citep{jung2020continual} that a significant fraction of ReLU units can become inactive across training and stay inactive under small perturbations. DNN can thus achieve natural sparsity after training, which can be leveraged to reduce the number of variables in the certification model \textit{cf.} right part of~\cref{fig:outline}. In modern architectures commonly used in computer vision, more than half of the neurons can be inactive~\citep{tjeng_evaluating_2019}, directly translating into an equivalent reduction in the number of variables. Beyond inactive units, stable \textbf{active} neurons should also be taken into account (see ~\citep{serra_lossless_2020, botoeva_efficient_2020}), thus achieving a far more compact formulation of the certification problem than in the original optimization problem.

The variable $z_k^{j}$ corresponding to a stable active neuron can be replaced by its linear expression $W_k^{j} z_{k-1} + b_k^{j}$. 
By recursively substituting each stable active neuron of vector $z_{k-1}$ by its linear expression in terms of $z_{k-2}$, and so on, the neuron $z_k^{j}$ can be linearly expressed across multiple layers (more information on Appendix~\ref{sup:stable_ablation}).

The chordal decomposition that we consider in ($SDP_T^j$) and ($SDP_U$) uses matrices $P_k$  which only model the links between two consecutive layers. 
However, when substituting the stable active neurons, unmodeled quadratic terms  may appear in Constraint~(\ref{ReLURef3}), \textit{i.e.} $ z^j_{k+1} (z^j_{k+1} - W_k z_k + b_k) = 0$. Indeed, let $z_{k}^{a}$ be the sub-vector of active neurons of layer $k$, 
and $z_{k}^{u}$ the sub-vector of unstable neurons (inactive neurons have been removed), Constraint~(\ref{ReLURef3}) becomes $z_{k+1}^{j} (z_{k+1}^{j} - W_{k,u}^{j} z_k^{u} - b_k) = z_{k+1}^{j}  W_{k,a}^{j} z_k^{a} = \sum\limits_{l=0}^{k-1} A_l^{i} z_l^{i} z_{k+1}^{j} + B z_{k+1}^{j}$, where $A$ and $B$ are derived from products of linear layer weights.  This formulation leads to new cross-layer dependencies as products $z_{k+1}^{j} z_l^{i}$ between non consecutive layers $l$ and $k+1$ appear, which are not represented in the $P_k$ matrices. To keep the chordal decomposition in our SDP relaxation, we bound the obtained right-hand side. As illustrated in the right part of~\cref{fig:outline}, we use McCormick cuts based on the bounds of $z_{k+1}^{j}$ and $z_l^{i}$ to create two linear upper bounds and two linear lower bounds on the quadratic products $z_{k+1}^{j} z_l^{i}$  :  \\
$\max \{LO^{1}(z_l^{i}, z_{k+1}^{j}), LO^{2}(z_l^{i}, z_{k+1}^{j})\} \leq z_{k+1}^{j} z_{l}^{i} \leq \min \{UP^{1}(z_l^{i}, z_{k+1}^{j}), UP^{2}(z_l^{i}, z_{k+1}^{j})\}$ 

where : 
\begin{equation*}
\left \{ \begin{array}{l}
    LO^{1}(z_l^{i}, z_{k+1}^{j})   =  U_{k+1}^{j} z_l^{i} + L_l^i z_{k+1}^{j} - U_{k+1}^{j}L_l^i  \quad   \nonumber \\
    LO^{2}(z_l^{i}, z_{k+1}^{j})  =   U_l^{i} z_{k+1}^{j}   \quad    \nonumber   \\
    UP^{1}(z_l^{i}, z_{k+1}^{j})   =   U_l^{i} z_{k+1}^{j} + U_{k+1}^{j} z_l^{i} - U_{k+1}^{j} U_{l}^{i} \quad  \nonumber \\
    UP^{2}(z_l^{i}, z_{k+1}^{j}) =  L_l^{i} z_{k+1}^{j} \quad \nonumber
\end{array} \right .
\end{equation*}

We then substitute the quadratic terms  between non consecutive layers by their linear bounds in the ReLU  constraint and obtain the following four relaxed ReLU  constraints for all $k+1 \in \mathcal{K}$, $j \in \mathcal{J}_{k+1}$:
  \begin{numcases}{}
   z_{k+1}^{j} (z_{k+1}^{j} - \!\!\!\!\!\!  \displaystyle{\sum_{\substack{u=1\\u \text{ unstable}}}^{n_k}} (W_{k+1,u}^{j} z_k^{u})) - b_{k+1}^j)   - B{\scriptstyle (j,k+1)} z_{k+1}^{j}\leq UP^r_{k+1,j}, \,\, r=1,2 \label{RelaxReluUP}\\
    z_{k+1}^{j} (z_{k+1}^{j} - \!\!\!\!\!\! \displaystyle{\sum_{\substack{u=1\\u \text{ unstable}}}^{n_k}} (W_{k+1,u}^{j} z_k^{u})) - b_{k+1}^j)  -B{\scriptstyle (j,k+1)}  z_{k+1}^{j}\geq LO^r_{k+1,j}, \,\, r=1,2 \label{RelaxReluLO}
  \end{numcases}
  
\noindent where $LO_1$, $LO_2$, $UP_1$, $UP_2$ are linear combinations of $z_{k+1}^{j}$ 
and of the unstable neurons across all layers $l < k$. 
More precisely, the upper and lower bounds in Constraints~\eqref{RelaxReluUP} and~\eqref{RelaxReluLO} are obtained as follows: \\

 $\left \{ \begin{array}{l}   LO^1_{k+1,j} = \sum\limits_{l=0}^{k-2} \sum\limits_{\substack{i \in \setJ, i \text{unstable}\\ A_l^{i} > 0}} LO^{1}(z_l^{i}, z_{k+1}^{j})
    + \sum\limits_{l=0}^{k-2} \sum\limits_{\substack{u \in \setJ, i \text{unstable}\\ A_l^{i} \leq 0}} UP^{1}(z_l^{i}, z_{k+1}^{j}) \nonumber \\
    LO^2_{k+1,j}= \sum\limits_{l=0}^{k-2} \sum\limits_{\substack{i \in \setJ, i \text{unstable}\\ A_l^{i} > 0}} LO^{2}(z_l^{i}, z_{k+1}^{j})
    + \sum\limits_{l=0}^{k-2} \sum\limits_{\substack{u \in \setJ, i \text{unstable}\\ A_l^{i} \leq 0}} UP^{2}(z_l^{i}, z_{k+1}^{j})  \nonumber 
\end{array} \right .$

$\left \{ \begin{array}{l}
    UP^1_{k+1,j} = \sum\limits_{l=0}^{k-2} \sum\limits_{\substack{i \in \setJ, i \text{unstable}\\ A_l^{i} > 0}} UP^{1}(z_l^{i}, z_{k+1}^{j})
    + \sum\limits_{l=0}^{k-2} \sum\limits_{\substack{u \in \setJ, i \text{unstable}\\ A_l^{i} \leq 0}} LO^{1}(z_l^{i}, z_{k+1}^{j}) \nonumber
   \\
   UP^2_{k+1,j} = \sum\limits_{l=0}^{k-2} \sum\limits_{\substack{i \in \setJ, i \text{unstable}\\ A_l^{i}> 0}} UP^{2}(z_l^{i}, z_{k+1}^{j})
    + \sum\limits_{l=0}^{k-2} \sum\limits_{\substack{u \in \setJ, i \text{unstable}\\ A_l^{i} \leq 0}} LO^{2}(z_l^{i}, z_{k+1}^{j})  \nonumber
\end{array} \right .$\\

 In our algorithm, we choose to only generate 4 
constraints via a deterministic selection, but in Constraints~\eqref{RelaxReluUP} and ~\eqref{RelaxReluLO}, depending on the sign of the coefficient $A_l^i$, each product can be relaxed by one of the two lower (respectively upper) bounds, resulting in a very large number of potential constraints (up to  $2^{n_0 + \sum_{l=1}^{k-1} n_l^{u}}$).

We prune active neurons on all layers except the penultimate layer, which is in the objective function of ($SDP_T^{j}$).
Finally, we obtain the following enhanced SDP relaxation, where $P_k$ matrices have been truncated: 
  \begin{numcases}{(SDP_U)}
    \min \,W_{K}^{y} P_{K-2}[z_{K-1}] + b_{K}^{y} -  \sum\limits_{j \in \setJKbar}  (W_{K}^{j} P_{K-2}[\beta_j z_{K-1}] + b_{K}^{j})   & \nonumber \\
    \mbox{s.t. } (\ref{ReLURef1_SDP})-(\ref{stab_SDP}),~ (\ref{SDP}),~(\ref{beta_sum1_SDP})-(\ref{RelaxReluLO}) \nonumber
\end{numcases} 
 Note that this pruning strategy of stable active neurons is a generic approach that can be applied to any SDP relaxation, either targeted or multiclass. As shown by our experiments of Section~\ref{exp:pruning}, applying this strategy to 
 ($SDP_U$) or ($SDP_T$) clearly speeds up the resolution. This size reduction comes at the cost of relaxing some equality constraints with inequalities in~\crefrange{RelaxReluUP}{RelaxReluLO}. However, the new constraints added to our formulation counterbalance this relaxation, ensuring that the overall certification performances remain competitive.

 Denoting by $n_k^{a}$ the number of stable active neurons on layer $k$, and $n_k^u$ the number of unstable neurons, Proposition~\ref{prop:size} specifies the reduction of size resulting from pruning.

\begin{proposition}
    The pruning of active neurons reduces the dimensions of each matrix variable $P_k$ for $k \in [K-3]$ from $(1 + n_k^{a} + n_k^{u} + n_{k+1}^{a} + n_{k+1}^{u})^{2}$ to $(1 + n_k^{u} + n_{k+1}^{u})^{2}$, and reduces $P_{K-2}$ in $SDP_U$ from $\big(1 + n_{K-2}^{a} + n_{K-2}^{u} + n_{K-1}^{a} + n_{K-1}^{u})^{2} + |\setJKbar| \big)^{2}$ to $\big(1 + n_{K-2}^{a} + n_{K-2}^{u} + n_{K-1}^{u})^{2} + |\setJKbar| \big)^{2}$.
    \label{prop:size}
\end{proposition}

    An interesting extension of our work would be to embed model $(SDP_U)$ within a Branch and Bound (B$\&$B) framework for complete certification. Indeed, classical B$\&$Bs handle one target class at a time, combining this algorithm with multiclass certification would be a first novelty. Moreover, most of B$\&$Bs are based on bound propagations and while fast to evaluate, the obtained lower bound is often weak, leading to excessive branching and tree deepening. A last advantage is that since deeper nodes exhibit more neuron stability, applying our pruning method will decrease the size $(SDP_U)$ making it cheaper and cheaper to solve. Note finally that solving faster SDP optimization problems is an active area of research (see~\cite{han_accelerating_2024}), what could further accelerate the whole  B$\&$B framework. \\

Note that our stable active pruning method can be generalized to other activation functions. Indeed, when the neuron lies entirely in its linear regime and can therefore be replaced by a linear decomposition on previous layer's outputs the method applies. This is the case for instance of the PReLU or of the Leaky ReLU functions, and more generally for any activation function that can be expressed as $\sigma(x)=\alpha x + \Delta,~\forall x \in [T,\infty)$ or $\sigma(x)=\alpha x + \Delta,~\forall x \in (-\infty,T]$.

\section{Experimental evaluation}

\subsection{Implementation Details}

We ran our experiments on a Linux machine on a 64-core CPU and 264 GB of RAM. We use the Python API of the MOSEK optimizer~\citep{mosek}. The bounds on the preactivation values are computed with $\alpha-\beta$-CROWN~\citep{wang_beta-crown_2021}  (see Appendices~\ref{supp_seq:algos} and ~\ref{supp_seq:exp} for details). We evaluate \texttt{SDP$_U$} (see details in Appendix~\ref{sup_seq:complements}) on \texttt{MNIST} (10 classes) and \texttt{CIFAR100} (100 classes). We reproduced the evaluation protocol from previous works on DNN certification~\citep{raghunathan_semidenite_2018,batten_efficient_2021,lan_tight_2022,li_sok_2023} by considering 3 different fully connected neural networks adversarially trained with \texttt{PGD} attacks (see Appendix~\cref{supp_seq:app_createDNN}). Neural networks used are 
  \texttt{9x100}~\citep{singh_abstract_2019} tested under the same $\epsilon = 0.026$;  \texttt{9x200}~\citep{singh_abstract_2019} tested under the same $\epsilon = 0.015$; and \texttt{FCNNA} on \texttt{CIFAR100}~\citep{li_sok_2023} tested under $\epsilon = 10^{-3}$. We have reproduced these networks to the best of our knowledge and report the detailed architecture and adversarial training in Appendix~\ref{supp_seq:exp}. We conducted our first two experiments on 100 data points: the first 10 images of each class from the \texttt{MNIST} train set. For \texttt{CIFAR100}, we experimented on 1000 data points using the same procedure: we selected the first 10 images of each class. The code is available at \url{\texttt{https://github.com/MargotBoyer/FastSDPCertification.git}}.

  \begin{table}[H]
\caption{Comparison of \texttt{SDP$_U$} with methods \texttt{SDP$_T$}, \texttt{SDP$_{T,layer}$}, \texttt{SDP$_T$}-IP.\\ 
$^\dagger$: mean time on 30 samples, $^\ddagger$: time on a single sample.}
\centering
\begin{small}
\begin{tabular}{@{}l@{~~~~}>{\columncolor{lblue}}c@{~~~~}>{\columncolor{lblue}}c@{~~~~}c@{~~~~}c@{~~~~}c@{~~~~}c@{~~~~}c@{~~~~}c@{}
}
    \hline
    \multirow{2}{*}{\textbf{Net.}}  & \multicolumn{2}{c}{\texttt{SDP$_U$}} (ours) & 
    \multicolumn{2}{c}{\texttt{SDP$_T$}} & 
    \multicolumn{2}{c}{\texttt{SDP$_{T,layer}$}} & 
    \multicolumn{2}{c}{\texttt{SDP$_T$}-IP}
     \\
     \cmidrule(r){2-3}\cmidrule(r){4-5}\cmidrule(r){6-7}\cmidrule(r){8-9}
     &  
     \textit{Cert.} &\textit{Time} & 
     \textit{Cert.} & \textit{Time}& 
     \textit{Cert.} & \textit{Time} &
     \textit{Cert.} & \textit{Time} \\

     \texttt{9x100} &  25/96
     & 5896 & 41/96 & 3699 & 41/96 & 559 & - & 6717$^\dagger$ \\

     \texttt{9x200}  & 71/100
     & 5054 & 56/100 & 12 184 & 56/100 & 906 & -  & 48 817$^\dagger$\\

     \texttt{FCNNA} &  121/139
     & 2091 & - & 92 776$^\ddagger$  & - & 72 136$^\ddagger$ & - & - \\

    \hline
    \hline
\end{tabular}
\end{small}

\label{tab:experiment1}
\end{table}

\subsection{State-of-the-art comparison}

To evaluate our method \texttt{SDP$_U$} with respect to other incomplete verifiers, we compare with the following methods: \texttt{SDP$_T$-IP}~\citep{raghunathan_semidenite_2018}, \texttt{SDP$_{T,layer}$}~\citep{batten_efficient_2021} and \texttt{SDP$_T$}~\citep{lan_tight_2022} use a chordal decomposition of matrices, ablation of stable inactive neurons, and the triangular constraint.
\texttt{SDP$_T$} additionally uses 10\% RLT cuts. For
\texttt{SDP$_U$} we use a chordal decomposition of matrices, ablation of stable inactive and active neurons, the triangular constraint, and $100\%$ of the RLT cuts for \texttt{9x100} and \texttt{FCNNA}, $60\%$ for \texttt{9x200}. 
The results are reported in Table~\ref{tab:experiment1}. Column \textit{Cert.} is the verified accuracy across \textit{all} targets: $x/y$ with $x$ the number of data formally certified as robust and $y$ the number of correctly classified data by the considered method. Column \textit{Time} is the mean total runtime per image (in seconds) across all classes.

We observe that \texttt{SDP$_U$} is significantly faster than methods \texttt{SDP$_T$} and \texttt{SDP$_T$}-IP on \texttt{9x200} and \texttt{FCNNA} networks. Due to the prohibitive CPU time required by \texttt{SDP$_T$}-IP, we ran the experiments on only 30 samples for \texttt{9x100} and  \texttt{9x200} networks. Moreover, for the same reason, we only ran on a single sample for the \texttt{FCNNA} network for methods \texttt{SDP$_T$} and \texttt{SDP$_{T,layer}$} (the solver failed for \texttt{SDP$_T$}-IP on this instance). Method \texttt{SDP$_{T,layer}$} is faster than our new method for networks with a small number of classes, but its CPU time grows exponentially with the number of classes. In contrast, \texttt{SDP$_U$} exhibits stable CPU time that is not affected by the number of classes. Regarding the certification rate, we observe that our new approach achieves the best results for 2 of the 3 networks considered. Specifically,  \texttt{SDP$_U$} outperforms  \texttt{SDP$_T$}  on the  \texttt{9x200} and the \texttt{FCNNA} networks,  and \texttt{SDP$_T$} does not scale to these two networks.

For the \texttt{9x200} network, this is partly due to the very large number of RLT constraints associated with the model, which cannot be fully incorporated. Although the use of RLT cuts in \texttt{SDP$_T$} improves the certification rate, it significantly slows down computation. Moreover, since \texttt{SDP$_T$} needs to solve up to 9 SDP models on \texttt{MNIST} for certification, adding RLT cuts becomes costly, and only a small proportion of them (10\% in our experiments) can be included while maintaining tractability. For \texttt{SDP$_T$}-IP, the absence of chordal decomposition leads to a prohibitive computational burden and certification could not be achieved in a reasonable time budget. 

Moreover, \texttt{SDP$_{T,layer}$} (as it is the case for \texttt{SDP$_T$}) does not scale when the number of classes becomes large, as observed on \texttt{FCNNA}. Finally, \texttt{SDP$_U$}  provides the best trade-off between certification performance and computation time. The aggregation of classes allows more RLT cuts to be added while remaining tractable. Additionally, pruning stable active neurons has a substantial impact on the size of the $(SDP_U)$ model, further reducing overall computation time. Consequently, our new method is particularly advantageous for larger networks or when dealing with a high number of classes compared with existing SDP-based approaches.

\subsection{Impact of the ablation of neurons}\label{exp:pruning}

\begin{figure}[h]
    \centering
    \begin{subfigure}{0.45\linewidth}
        \centering
        \includegraphics[width=3.5cm]{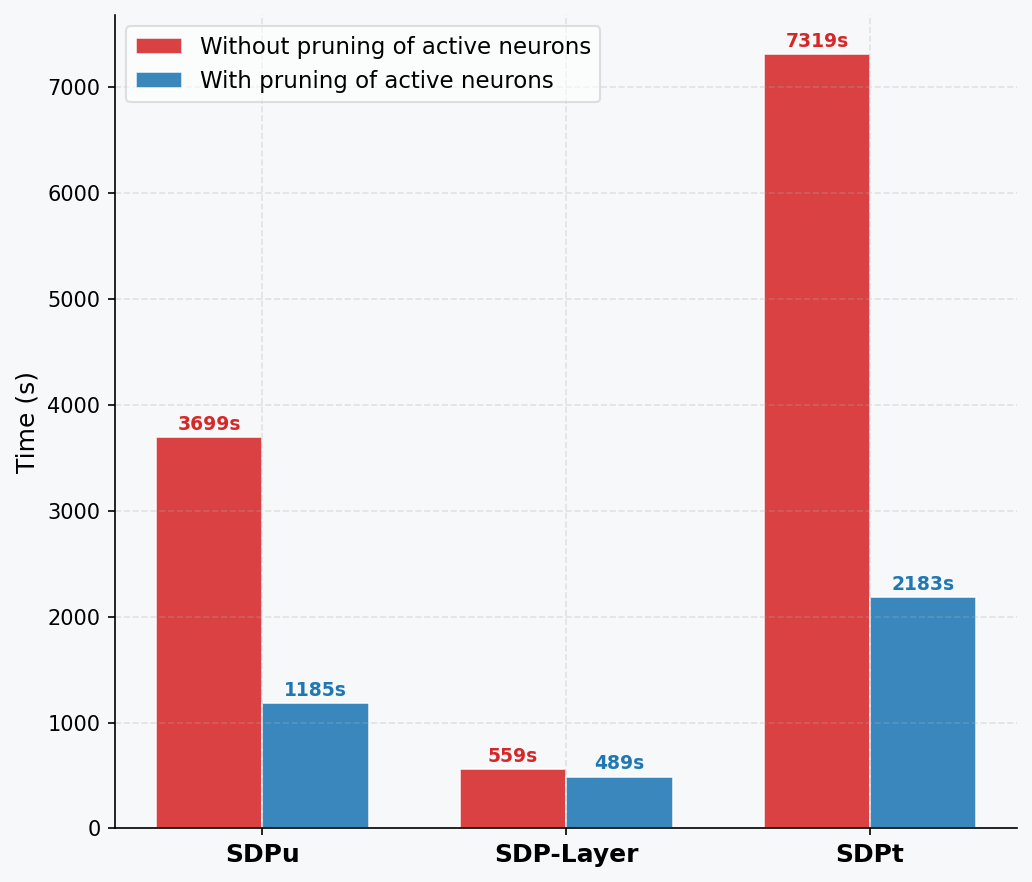}
        \caption{Execution Time.}
        \label{fig:image1}
    \end{subfigure}
    \begin{subfigure}{0.45\linewidth}
        \centering
        \includegraphics[width=3.5cm]{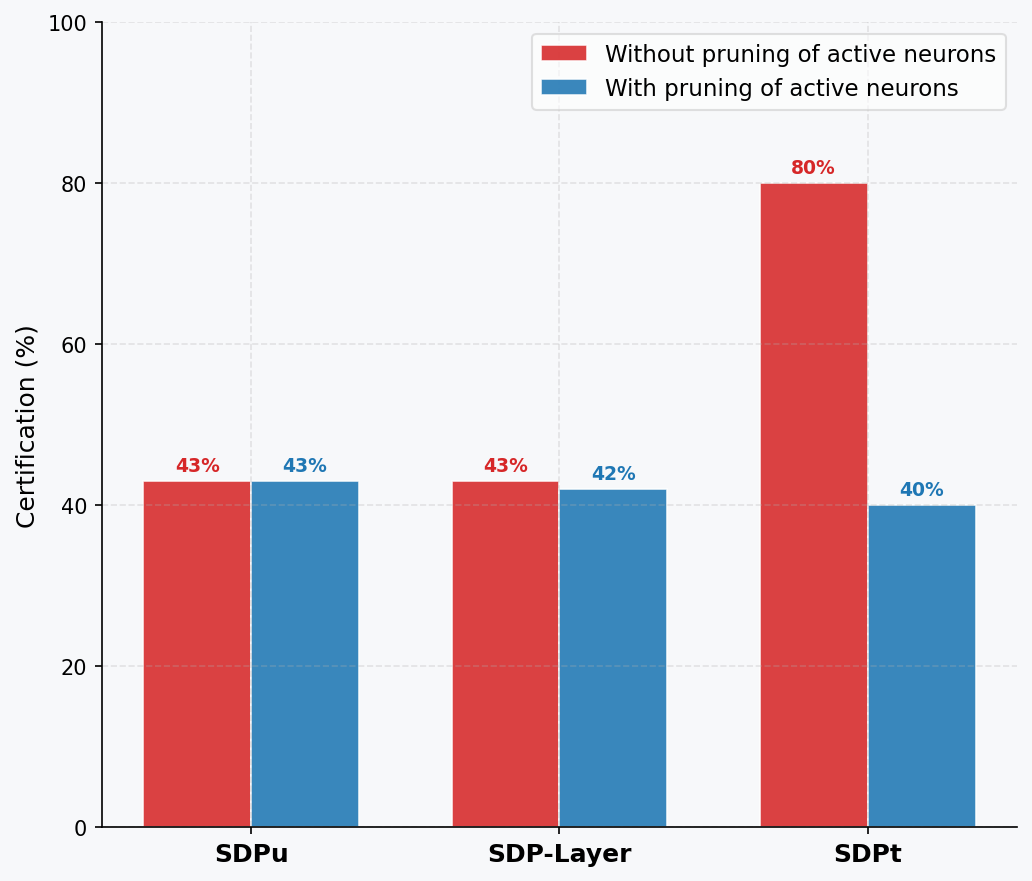}
        \caption{Certification Rate.}
        \label{fig:image2}
    \end{subfigure}\\

    \begin{subfigure}{0.47\linewidth}
        \centering
    \includegraphics[width=4.5cm]{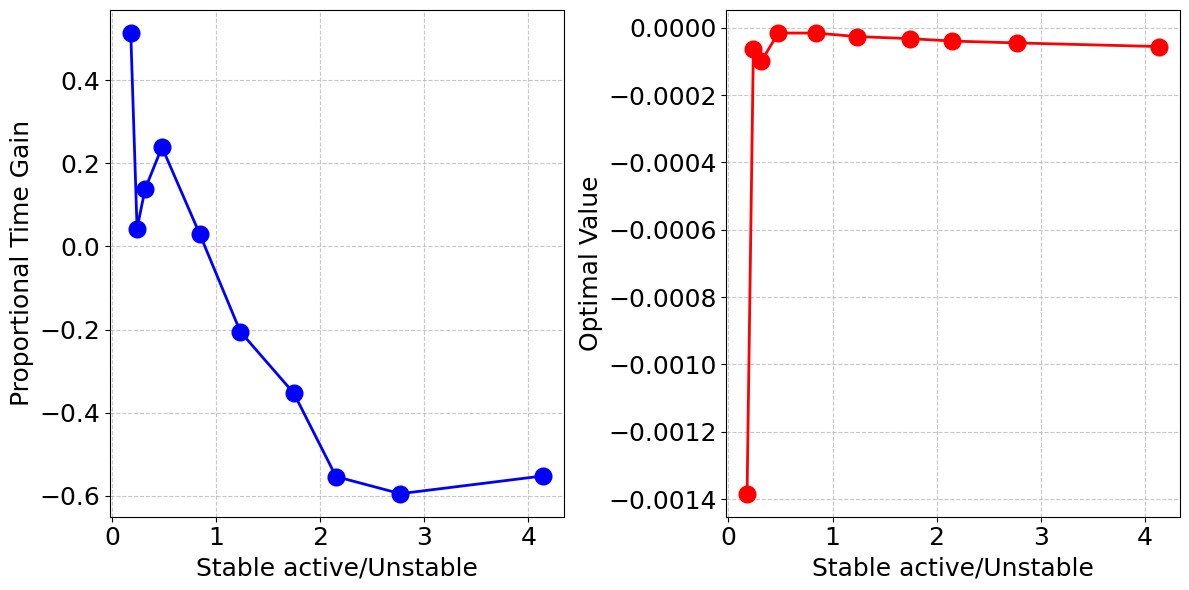}
        \caption{\texttt{SDP-Layer} study.}
    \label{fig:pruning_study_sdp_layer}
    \end{subfigure}
    \begin{subfigure}{0.47\linewidth}
        \centering
    \includegraphics[width=4.5cm]{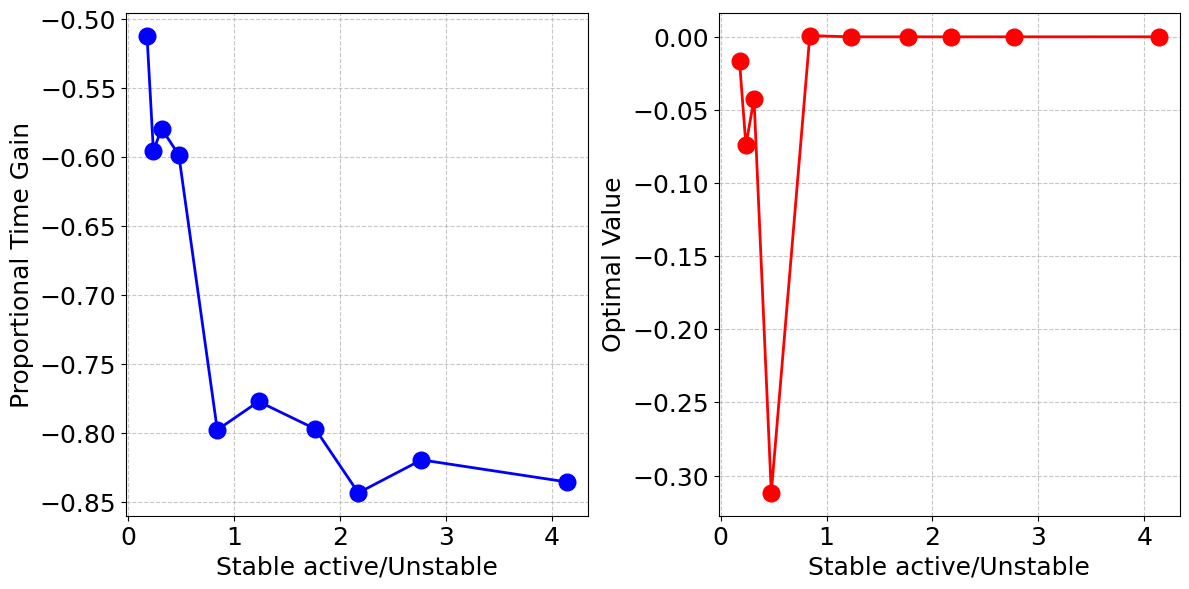}
        \caption{\texttt{SDP$_T$} study.}
        \label{fig:pruning_study_sdpt}
    \end{subfigure}\\

    \begin{subfigure}{0.47\linewidth}
        \centering
    \includegraphics[width=4.5cm]{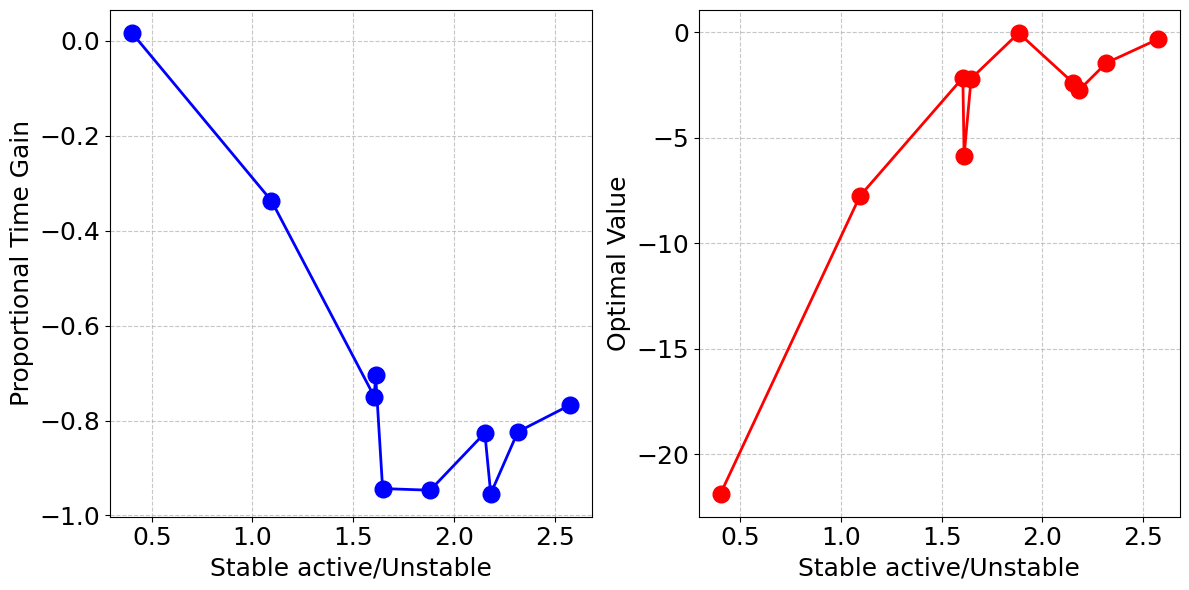}
        \caption{\texttt{SDP$_U$} study.}
        \label{fig:pruning_study_sdpu}
    \end{subfigure}

    \caption{Impact of stable active neuron pruning on SDP-based certification. 
    (a-b)~CPU time and certification rates. (c-d-e)~Proportional time gain and 
    objective value decrease as a function of the ratio of stable active to unstable neurons.}
    \label{fig:both}
\end{figure}

We now study the impact of the ablation of {stable} active neurons on the performance of methods \texttt{SDP$_T$}, \texttt{SDP-Layer} and \texttt{SDP$_U$}. Since the CPU time without pruning is large for \texttt{SDP$_U$}, we ran our experiments on network \texttt{9x100} on $10$ randomly picked data for \texttt{SDP$_U$}. We consider two configurations: pruning of inactive neurons, and pruning of both stable active and inactive neurons.

We start by evaluating the impact of pruning active neurons on CPU time and percentage of certification. We plot the results in Figures~\ref{fig:image1} and~\ref{fig:image2}. In terms of CPU time, the results reveal a similar trend for the three methods. As expected, the greater the number of pruned neurons, the faster the resolution. Indeed, performing this ablation along with chordal decomposition allows for drastically reducing the number of variables and constraints. For this small network, the percentage of certification is more impacted with \texttt{SDP$_U$} since the pruning relaxes some of the constraints. Note that this degradation could be mitigated with tighter bounds on stable active neurons, and when applied on larger networks where the percentage of \texttt{$SDP_U$} with pruning is the highest within the 3 compared methods. 

We now examine how the percentage of neurons removed by our pruning method influences the performances. Figures~\ref{fig:pruning_study_sdp_layer}, \ref{fig:pruning_study_sdpt}, and \ref{fig:pruning_study_sdpu} present the proportional CPU time gain (blue) and the optimal value decrease (red) when using our pruning as a function of the ratio of active to unstable neurons. 
Each point on a curve corresponds to a given ratio of active to unstable neurons and reports the associated proportional CPU time gain (blue) or decrease in the optimal value (red). When the problem has greater neuron stability (i.e. when stable active neurons outnumber unstable ones) our pruning significantly reduces computation time (by more than 50\% on \texttt{SDP$_T$} for any stability ratio, by more than 40\% on \texttt{SDP-Layer} for a ratio above 2, and by more than 60\% on \texttt{SDP$_U$} for a ratio greater than 1.5). Moreover, the loss in optimal value is significantly reduced as neuronal stability increases.

\subsection{Impact of the aggregation of classes}\label{exp:aggregation}

 Finally, we further assess the scalability of our method with respect to large-scale, multi-class datasets. We construct a composite dataset by merging \texttt{EMNIST Balanced}, \texttt{KMNIST}, and \texttt{FashionMNIST}, resulting in a total of 67 distinct classes. We train neural networks on subsets of this dataset, with 5, 20, 50, and 67 classes respectively, with one representative for each class and $\epsilon=0.01$. 
 
\begin{figure}[H]
    \centering
\includegraphics[width=3.1cm]{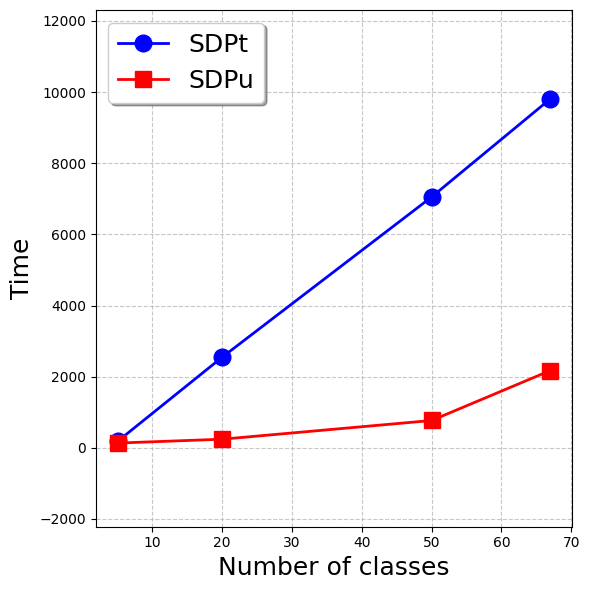}
        \caption{Time vs \# of classes}
        \label{fig:image3}
\end{figure}

We compare the runtime performance of \texttt{SDP$_U$} and \texttt{SDP$_T$} in Figure~\ref{fig:image3}, where each line plots the computation time (in seconds) according to the number of classes. We observe that the computation time of \texttt{SDP$_T$} increases greatly with the number of classes, while the computation time of \texttt{SDP$_U$} remains stable. 
Note that certification itself becomes intrinsically harder as the number of classes grows, since by Theorem~\ref{prop:equi} the certified bound is a minimum taken over an increasing number of per-class quantities. 
This is a property of the certification problem itself, not a limitation specific to \texttt{SDP$_U$}.
Moreover, on the computational side, it is well known that solving large SDP models can be time consuming. 
However, SDP solvers are improving fast, with promising GPU-accelerated implementations~\cite{han_accelerating_2024} now emerging, clearly encouraging the aggregation of classes for large multi-class datasets.

\section{Conclusion}

We introduced a new multiclass SDP model to certify ReLU networks across all targets in a single model. Our model enables a significant speedup compared to single targeted SDP models, in particular for datasets with a large number of classes.  We are further able to improve both targeted and untargeted models thanks to a size reduction of the SDP models by pruning variables corresponding to stable active neurons. A first perspective will be to extend this pruning method to other activation functions. Moreover, since the bottleneck of SDP-based methods is its computation time on large networks,  a relevant perspective will be to accelerate the solution of $(SDP_U)$ using a conic bundle method.

\begin{credits}
\subsubsection{\ackname} 
This work was performed using HPC resources from GENCI–IDRIS (Grant 2025-AD011016715)
\subsubsection{\discintname}
 The authors declare that they have no competing interests.
\end{credits}
%
%

\bibliographystyle{splncs04}
\bibliography{Certification_Problem_bibtex}

\nocite{Cor76, ehlers_formal_2017, katz_reluplex_2017, gurfinkel_marabou_2024, fischetti_deep_2018, dsouza_maximum_2017, liao_solution-aware_2025}
\nocite{bunel_branch_2020,ferrari_complete_2022,jaeckle_neural_2021,lu_neural_2019, henriksen_deepsplit_2021, zhou_scalable_2024, zhang_general_2022, lan_semidefinite_2023,chiu_sdp-crown_2025}
\nocite{madry_towards_2019, goodfellow_explaining_2015, kurakin_adversarial_2017, papernot_limitations_2015, moosavi-dezfooli_deepfool_2016, carlini_towards_2017, salman_provably_2019}

\end{document}


\title{Fast SDP certification of neural networks : towards large multi-class datasets}

\titlerunning{Fast SDP certification of neural networks : towards large multi-class datasets}

\authorrunning{M. Boyer et al.}




\appendix
\onecolumn

\section{Proofs}
\label{supp_seq:proofs}
\subsection{Proof of Proposition 3}
\begin{proof}
Constraint~(\ref{boundRef}) taking $k=0$ is equivalent to Constraint~(\ref{ball}), ensuring that $z_0\in\mathcal{B}_{\epsilon}(x)$. Thus, $v(QP_{T}^j) \geq 0$ implies $\min_{z_0\in\mathcal B_\epsilon(x)} z_K^{y}-z^j_K\geq 0$.
\end{proof}

\subsection{Proof of Theorem 1}
\begin{proof}
\begin{enumerate}
\item  We first prove that $v(QP_{T}^{\bar j})\geq v(QP_{U})$ for all $\bar j\in\setJKbar$. Let $z^*$ be an optimal solution of $(QP_{T}^{\bar j})$. We build the solution $(z=z^*, \beta)$ of $(QP_{U})$ with $\beta_{\bar j} = 1$ and $\beta_j = 0 \; \;  \forall j\in\setJKbar\backslash\{\bar j\}$ satisfying Constraints~(\ref{binary}) and (\ref{beta_sum1}). Obviously, the two objective functions have the same value.
\item Then, we prove that $ \min\limits_{\bar j \in \setJKbar} v(QP_{T}^{\bar j})\leq v(QP_{U})$. \\
Let $(z^*, \beta^*)$ be an optimal solution of $(QP_{U})$ with $\beta_{\bar j}=1$ for $\bar j\in\setJKbar$. Obviously, $z=z^*$ is feasible for $(QP_{T}^{\bar j})$.  Once again, the objective values are the same.
\end{enumerate}
\end{proof}

\subsection{Proof of Proposition 5}
\begin{proof}
    We ensure that the inequalities are valid in all three possible cases: 
    \begin{enumerate}
        \item If $\beta_{j_1} = 1$ and  $\beta_{j_2} = 0$, then both sides of~(\ref{eq:cut1}) are equal to $0$ and~(\ref{eq:cut2}) becomes $U^{j_1}\geq z_K^{j_1}(x)$.        
        \item $\beta_{j_1} = 0$ and $\beta_{j_2} = 1$, (\ref{eq:cut1}) and~(\ref{eq:cut2}) respectively lead to $z_K^{j_2}(x)\leq U_{K}^{j_2}$ and $z_K^{j_2}(x)\geq z_K^{j_1}(x)$. The latter is valid since when $\beta_{j_2}=1$, class $j_2$ provides the worst example.

        \item Otherwise, $\beta_{j_1} = \beta_{j_2} = 0$ and we obtain $L_{K}^{j_1}\leq z_K^{j_1}(x)$ from~(\ref{eq:cut1}) and $U_K^{j_1}\geq z_K^{j_1}(x)$ from~(\ref{eq:cut2}).
    \end{enumerate}
\end{proof}

\subsection{Proof of Proposition 6}

\begin{proof}
    
        We consider here the case where the neural network has more than one hidden layer, \textit{i.e.}, $K > 2$:
        \begin{itemize}
            
            \item If $k \in [K-3]$ for $SDP_U$, or if $k \in [K-2]$ for $SDP_T$: \\
            The matrix $P_k$ in $SDP_T$ and $SDP_U$ is of size $1 + n_k^{u} + n_k^{a} + n_{k+1}^{u} + n_{k+1}^{a}$, where $n_k^{a}$ is the number of stable active neurons in layer $k$, and $n_k^{u}$ the number of unstable neurons. Note that $n_k^{a} + n_k^{u}$ corresponds to the layer size $n_k$ minus the number of stable inactive neurons. Pruning the stable active neurons in layer $k$ removes $n_k^{a}$ rows from the matrix $P_k$. The same applies to layer $k+1$. The final matrix $P_k$ is therefore of size $1 + n_k^{u} + n_{k+1}^{u}$. The number of matrix entries thus decreases from $(1 + n_k^{u} + n_k^{a} + n_{k+1}^{u} + n_{k+1}^{a})^{2}$ to $(1 + n_k^{u} + n_{k+1}^{u})^{2}$.
            \item If $k = K-2$ for $SDP_U$: the matrix size without pruning of stable active neurons is $1 + n_{K-2}^{u} + n_{K-2}^{a} + n_{K-1}^{u} + n_{K-1}^{a} + |\setJKbar|$. Removing the stable active neurons reduces the matrix size to $1 + n_{K-2}^{u} + n_{K-1}^{u} + |\setJKbar|$. The number of matrix entries thus goes from $(1 + n_{K-2}^{u} + n_{K-2}^{a} + n_{K-1}^{u} + n_{K-1}^{a} + |\setJKbar|)^{2}$ to $(1 + n_{K-2}^{u} + n_{K-1}^{u} + |\setJKbar|)^{2}$.
            
        \end{itemize}
        
        If the neural network has one or zero hidden layer, there is only one matrix variable in the model, namely $P_0$. The same reasoning as above can be applied to compute the reduction in the number of matrix entries.
        
\end{proof}

\subsection{Constraints~\cref{boundRef_SDP_k} or~\cref{boundRef_SDP} also implies $L_k \leq P[z_k] \leq U_k$ in the SDP case}
\begin{proof}
We want to prove that $P[z_k z_k^T] - (U_k + L_k) P[z_k] + U_k L_k \leq 0$ implies $L_k \leq z_k \leq U_k$. \\
Assume $P$ is positive semi-definite, so every principal minor (ie. every determinant of any principal submatrix) is nonnegative. In particular :

$det\begin{pmatrix} P[1] & P[z_k] \\
    P[z_k] & P[z_k z_k^T]
\end{pmatrix} \geq 0$, and since $P[1] = 1$, this gives $P[z_k]^{2} \leq P[z_k z_k^T]$. Plugging $P[z_k]^{2}$ into~\cref{boundRef_SDP_k} we obtain $P[z_k]^{2} - (U_k + L_k) P[z_k] + U_k L_k \leq 0$, which is equivalent to $(U_k - P[z_k]) (P[z_k] - L_k) \geq 0$. $L_k \leq U_k$ gives : $L_k \leq P[z_k] \leq U_k$.

The same reasoning applies with $P_k$.
\end{proof}
\newpage

\newpage
\section{Complements on method}
\label{sup_seq:complements}

\subsection{Algorithms comparison}
\label{supp_seq:algos}
\subsubsection{Full certification algorithm}

This subsection describes the procedure used to prune trivially certified targets (certified with the bounds on the logits obtained via $\alpha-\beta-$CROWN). To ensure a fair comparison, targeted SDPs are only run on targets that are not trivially certified.  Note that \texttt{9x100} and \texttt{9x200} exhibited nearly no trivial certifications, meaning almost all 9 targets were evaluated per input, whereas \texttt{FCNNA} saw effective pruning, cutting the number of evaluated targets from 99 down to 93 for the first sample.
\begin{paddedalgo}
\SetAlgoVlined
\caption{Full certification of $\epsilon$-robustness of a DNN with targeted SDP model $(SDP_T)$ on data $x \in \setX$}
\label{algo:fullCert}
\texttt{
SDP$_T$($\epsilon,x$)} \\
cert $\gets$ true\;
  \For{$j \in \setJKbar$}{
    Compute bounds with $\beta$-CROWN. \\
    \If{$\exists l \in \setJKbar\cup{y}, \; \; U_K^{j} < L_K^{l}$}{
      \textbf{break}\;}
    $v^*_j \gets$ Solve $(SDP_{T}^{\epsilon,j,x})$\;
    \If{$v^*_j \leq 0$}{
      cert $\gets$ false\;
    }
  }
  
\Return cert\;

\end{paddedalgo}

\begin{paddedalgo}
\SetAlgoVlined
\caption{Full Certification of $\epsilon$-robustness of a DNN with untargeted SDP model ($
SDP_U$) on data $x \in \setX$}
\label{algo:fullCert}
\texttt{SDP$_U$($\epsilon,x$}) \\
cert $\gets$ true\;
Compute bounds with $\beta$-CROWN. \\
T $\gets \{j \in \setJKbar, \;\;\; \exists l \in \setJKbar\cup{y}, \;\;\; U_K^{j} < L_K^{l}\} $ \\
$v^* \gets$ Solve $(SDP_{U}^{\epsilon,T,x})$\;
\If{$v^* \leq 0$}{
  cert $\gets$ false\;
}
\Return cert\;
\end{paddedalgo}

\subsection{Stable active neurons ablation}
\label{sup:stable_ablation}

\subsubsection{Presentation}

In this section, we explain with more details our ReLU  constraint relaxation in the context of stable active neurons ablation. 

For a given $k \in [K-2]$, the quadratic ReLU  constraint on neuron $j$ of layer $k+1$ is 

$$z_{k+1}^{j} \Big (z_{k+1}^{j} - \displaystyle{\sum_{i =1:~i \text{ unstable or active}}^{n_k} } (W_{k+1,i}^{j} z_k^{i}) -  b_{k+1}^j \Big) = 0$$
or equivalently, by separating unstable and stable active neurons:
$$z_{k+1}^{j} (z_{k+1}^{j} - \displaystyle{\sum_{u= 1:~u\text{ unstable}}^{n_k}} (W_{k+1,u}^{j} z_k^{u})) - b_{k+1}^j) = z_{k+1}^{j} ( \displaystyle{\sum_{a=1:~a \text{ active}}^{n_k}} W_{k+1,a}^{j} z_k^{a})$$

\noindent \textit{i) Linear expression of stable active neurons}\\

Note that each stable active neuron $a$ of layer $k$ can be formulated as a linear combinations of previous layers unstable outputs $u$, \emph{i.e.}
$$z_k^{a}~=~\sum\limits_{l=0}^{k-1}~\sum\limits_{u=1:~u \text{ unstable}}^{n_l}~\lambda_l^{u}{\scriptstyle (a,k)}~z_l^{u}~+~\gamma{\scriptstyle (a,k)}$$
where $\gamma{\scriptstyle (a,k)} \in \mathbb{R}$ can be viewed as the bias of the expression of $z_k^a$, and $\lambda_l^{u}{\scriptstyle (a,k)} \in\mathbb{R}$ 
is the linear coefficient of neuron $u$ of layer $l$ in the expression of $z_k^a$. We thus have to compute for each \emph{stable active} neuron $a$ of layer $k$, and for all \emph{unstable} neuron $u$ of layer $l \in [0,k-1]$ the coefficients $\lambda_l^{u}{\scriptstyle (a,k)}$ and the bias term $\gamma(a,k)$.

For each unstable neuron $u$ of layer $k$, these coefficients can then be computed dynamically as follows : 
\begin{itemize}
    \item if $k=1$, and since all neurons $u$ of layer $0$ are \emph{unstable}, we have for all  stable active neuron $a$ of layer $1$:  $$z_1^{a} = \displaystyle{\sum_{i = 0}^{n_0}} W_{1,i}^{a} z_0^i + b_1^{a}$$ and therefore $\lambda_0^{u}{\scriptstyle (a,1)}=W_{1,u}^{a}$, and $\gamma{\scriptstyle (a,1)}=b_1^{a}$;
    \item if $k>1$, for a stable active neuron $a$ on such layer, 
    \begin{small}
    \begin{align*}
        z_k^{a} &= \displaystyle{\sum_{i=1:~i \text{ stable or unstable}}^{n_{k-1}}}  W_{k,i}^{a} z_{k-1}^i + b_k^{a} \\
        &= \displaystyle{\sum_{u=1:~u \text{ unstable}}^{n_{k-1}}} W_{k,u}^{a} z_{k-1}^{u} + \displaystyle{\sum_{r=1:~r \text{ active }}^{n_{k-1}}} W_{k,r}^{a} z_{k-1}^{r} + b_k^{a} \\
        &= \displaystyle{\sum_{u=1:~u \text{ unstable}}^{n_{k-1}}} W_{k,u}^{a} z_{k-1}^{u} \\&\quad+ \displaystyle{\sum_{r=1:~r \text{ active }}^{n_{k-1}}}
         W_{k,r}^{a}  \Big(\sum\limits_{l=0}^{k-2} \sum\limits_{u=1:~u \text{ unstable}}^{n_l} 
         \lambda_{l}^{u}{\scriptstyle (r,k-1)} z_{l}^{u} + \gamma{\scriptstyle (r,k-1)} \Big) + b_k^{a}  \\
        &= \displaystyle{\sum_{u=1:~u \text{ unstable}}^{n_{k-1}}} \underbrace{W_{k,u}^{a}}_{\lambda^u_{k-1}(a,k)} z_{k-1}^{u} +       
           \sum\limits_{l=0}^{k-2} \sum\limits_{u=1:~u \text{ unstable}}^{n_l} \underbrace{\displaystyle{\sum_{r=1:~r\text{ active }}^{n_{k-1}}} W_{k,r}^{a}\lambda_{l}^{u}{\scriptstyle (r,k-1)}}_{\lambda_l^u(a,k)} z_{l}^{u} \\&\quad+ \underbrace{\displaystyle{\sum_{r=1:~r \text{ active }}^{n_{k-1}}} W_{k,r}^{a} \gamma{\scriptstyle (r,k-1)}  + b_k^{a}}_{\gamma(a,k)}
    \end{align*}
    \end{small}
    Which gives  : 
    For each stable active neuron $a$ of layer $k>1$:
        \begin{itemize}
         \item $\gamma{\scriptstyle (a,k)} = \displaystyle{\sum_{r=1:~r \text{ active }}^{n_{k-1}}} W_{k,r}^{a} \gamma{\scriptstyle (r,k-1)}  + b_k^{a}$
    \item For all unstable neuron $u$ on layer $k-1$ :
        $\lambda_{k-1}^{u}{\scriptstyle (a,k)} = W_{k,u}^{a}$ 
        \item For all unstable neuron $u$ on layer $l \leq k-2$ :
        
        $\lambda_{l}^{u}{\scriptstyle (a,k)} = \displaystyle{\sum_{r=1:~r \text{ active }}^{n_{k-1}}} W_{k,r}^{a}\lambda_{l}^{u}{\scriptstyle (r,k-1)}$ 
    
        \end{itemize}
     
\end{itemize}

\noindent \textit{ii) Rewriting of the ReLU quadratic constraint}\\

Using the above coefficients, our ReLU  quadratic constraint can be rewritten

    \begin{small}
    $\begin{array}{l}
z_{k+1}^{j} (z_{k+1}^{j} - \displaystyle{\sum_{u=1:~u\text{ unstable}}^{n_k}} (W_{k+1,u}^{j} z_k^{u})) - b_{k+1}^j) = z_{k+1}^{j} ( \displaystyle{\sum_{a=1:~a \text{ active}}^{n_k}} W_{k+1,a}^{j} z_k^{a}) \\
 \quad\quad\;=  \displaystyle{\sum_{a=1:~ a\text{ active}}^{n_k}} W_{k+1,a}^{j} z_{k+1}^{j} z_k^{a}  \\
\quad\quad\;=  \displaystyle{\sum_{a=1:~a \text{ active}}^{n_k}} W_{k+1,a}^{j} z_{k+1}^{j}   \Big(\sum\limits_{l=0}^{k-2} \sum\limits_{u=1:~u \text{ unstable}}^{n_l} 
         \lambda_{l}^{u}{\scriptstyle (a,k)} z_{l}^{u} + \gamma{\scriptstyle (a,k)} \Big)\\
\quad\quad\;=    \sum\limits_{l=0}^{k-2} \sum\limits_{u=1:~u \text{ unstable}}^{n_l}\underbrace{ \displaystyle{\sum_{a=1:~a \text{ active}}^{n_k}} W_{k+1,a}^{j} 
         \lambda_{l}^{u}{\scriptstyle (a,k)}}_{A_l^{u}{\scriptstyle (j,k+1)} } z_{k+1}^{j} z_{l}^{u} \\
         \quad\quad\;+ \underbrace{\displaystyle{\sum_{a=1:~a \text{ active}}^{n_k}} W_{k+1,a}^{j}  \gamma{\scriptstyle (a,k)}}_{B{\scriptstyle (j,k+1)} } z_{k+1}^{j}\\[5ex]
        \end{array}$  
    \end{small}

The ReLU quadratic constraint becomes : \\
\begin{small}
    $\begin{array}{l}   
       z_{k+1}^{j} (z_{k+1}^{j} - \!\!\!\!\!\!\!\!\!\!\!\!
       \displaystyle{\sum_{u=1:~u \text{ unstable}}^{n_k}} (W_{k+1,u}^{j} z_k^{u})) - b_{k+1}^j)  = \sum\limits_{l=0}^{k-2} \sum\limits_{u=1:~u \text{ unstable}}^{n_l} A_l^{u}{\scriptstyle (j,k+1)}  z_{k+1}^{j} z_{l}^{u} + B{\scriptstyle (j,k+1)}  z_{k+1}^{j} \nonumber
    \end{array}$ 
\end{small}

Or equivalently:
\begin{small}
\begin{equation}
 z_{k+1}^{j} (z_{k+1}^{j} - \displaystyle{\sum_{u=1:~u \text{ unstable}}^{n_k}} (W_{k+1,u}^{j} z_k^{u})) - b_{k+1}^j)  -B{\scriptstyle (j,k+1)}  z_{k+1}^{j} = \sum\limits_{l=0}^{k-2} \sum\limits_{u=1:~u \text{ unstable}}^{n_l} A_l^{u}{\scriptstyle (j,k+1)}  z_{k+1}^{j} z_{l}^{u} \label{ReluRelaxed}
\end{equation}
 \end{small}

where for each unstable neuron $j$ of layer $k+1$, we have:
  \begin{itemize}
\item $A_l^{u}{\scriptstyle (j,k+1)}  = \sum\limits_{a=1 \text{ active}}^{n_k} W_{k+1,a}^{j} \lambda_l^{u}{\scriptstyle (a,k)}$ for all $u$ unstable $\in \{1,...,n_l\}, l \leq k-1$
\item  $B{\scriptstyle (j,k+1)} =  \sum\limits_{a=1 \text{ active}}^{n_k} W_{k+1,a}^{j} \gamma{\scriptstyle (a,k)}$
  \end{itemize}

\noindent \textit{iii) Relaxation of the ReLU constraints}\\

We now consider the quadratic terms $z_{k+1}^{j} z_l^{u}$ involved in Constraints~(\ref{ReluRelaxed}). Observe that for a neuron $j$ of layer $k+1$, the variables $z_l^{u}$ stands for a layer $l$ that precedes layer $k$ by a distance of at least 2, \emph{i.e.}, $l \in [k-1]$, and $u$ is an unstable neuron of the associated layer $l$. Since we use chordal decomposition in our SDP relaxation,  only products of neurons of two consecutive layers are modeled and the quadratic terms $z_{k+1}^{j} z_l^{u}$ are not considered in our formulation. Consequently, we cannot directly add Constraint~(\ref{ReluRelaxed}) in our SDP. To keep a ReLU constraint in our formulation, our idea is therefore to upper and lower bound the left-hand side of Constraint~(\ref{ReluRelaxed}) using the McCormick envelopes~\citep{Cor76} of the products $z_{k+1}^{j} z_l^{u}$.

We have $0 \leq z_{k+1}^{j} \leq U_{k+1}^{j}$, $L_l^u \leq z_{l}^{u} \leq U_{l}^{u}$. Note that we compute  $U_{k+1}^{j}$, $U_{l}^{u}$ and $L_l^u$ with \texttt{$\beta$-CROWN}, and for $l\geq 1$, we take $L_l^u = 0$ as it is a better bound than the lower bound on the preactivation vector. The McCormick envelopes stand for:
    \begin{align}{}
    z_{k+1}^{j} z_{l}^{u}  &\leq  \min \{ U_{k+1}^{j} z_l^{u} + L_l^u z_{k+1}^{j} - U_{k+1}^{j}L_l^u, U_l^{u} z_{k+1}^{j} \}  \quad   \label{stab:bound1}\\
    z_{k+1}^{j} z_{l}^{u}  &\geq  \max \{L_l^{u} z_{k+1}^{j} , U_l^{u} z_{k+1}^{j} + U_{k+1}^{j} z_l^{u} - U_{k+1}^{j} U_{l}^{u} \} \quad \label{stab:bound2}
     \end{align}
     
To upper bound Constraints~(\ref{ReluRelaxed}), we generate the 4 following expressions:
    \begin{scriptsize}
    \begin{numcases}{}
    UP_1 =  \displaystyle{\sum_{l=0}^{k-2}} \displaystyle{\sum_{\substack{u=1:~u \text{ unstable} \\ A_l^{u}{\scriptstyle (j,k+1) >0}}}^{n_l}} A_l^{u}{\scriptstyle (j,k+1)} ( U_{k+1}^{j} z_l^{u} + L_l^u z_{k+1}^{j} - U_{k+1}^{j}L_l^u) \nonumber\\\quad\quad\;- \displaystyle{\sum_{l=0}^{k-2}} \displaystyle{\sum_{\substack{u=1:~u \text{ unstable} \\ A_l^{u}{\scriptstyle (j,k+1) <0}}}^{n_l}} A_{k+1}^{j}{\scriptstyle (j,k+1)}  L_l^{u} z_{k+1}^{j} \nonumber \\ 
    UP_2 =  \displaystyle{\sum_{l=0}^{k-2}} \displaystyle{\sum_{\substack{u=1:~u \text{ unstable} \\ A_l^{u}{\scriptstyle (j,k+1) >0}}}^{n_l}} A_{k+1}^{j}{\scriptstyle (j,k+1)}  U_l^{u} z_{k+1}^{j} \\\quad\quad\;- \displaystyle{\sum_{l=0}^{k-2}} \displaystyle{\sum_{\substack{u=1:~u \text{ unstable} \\ A_l^{u}{\scriptstyle (j,k+1) <0}}}^{n_l}} A_{k+1}^{j}{\scriptstyle (j,k+1)} (U_l^{u} z_{k+1}^{j} + U_{k+1}^{j} z_l^{u} - U_{k+1}^{j} U_{l}^{u}) \nonumber\\
     \end{numcases} 
    \end{scriptsize}

       To lower bound Constraints~(\ref{ReluRelaxed}), we generate the 4 following expressions:
       
    \begin{scriptsize}
    \begin{numcases}{}
    LO_1 =  \displaystyle{\sum_{l=0}^{k-2}} \displaystyle{\sum_{\substack{u=1:~u \text{ unstable} \\ A_l^{u}{\scriptstyle (j,k+1) > 0}}}^{n_l}} A_{k+1}^{j}{\scriptstyle (j,k+1)}  L_l^{u} z_{k+1}^{j} \nonumber\\\quad\quad\;- \displaystyle{\sum_{l=0}^{k-2}} \displaystyle{\sum_{\substack{u=1:~u \text{ unstable} \\ A_l^{u}{\scriptstyle (j,k+1) < 0}}}^{n_l}} A_l^{u}{\scriptstyle (j,k+1)} ( U_{k+1}^{j} z_l^{u} + L_l^u z_{k+1}^{j} - U_{k+1}^{j}L_l^u)  \nonumber \\ 
    LO_2 =  \displaystyle{\sum_{l=0}^{k-2}} \displaystyle{\sum_{\substack{u=1:~u \text{ unstable} \\ A_l^{u}{\scriptstyle (j,k+1) > 0}}}^{n_l}} A_{k+1}^{j}{\scriptstyle (j,k+1)}  U_l^{u} z_{k+1}^{j} \\\quad\quad\;- \displaystyle{\sum_{l=0}^{k-2}} \displaystyle{\sum_{\substack{u=1:~u \text{ unstable} \\ A_l^{u}{\scriptstyle (j,k+1) <0}}}^{n_l}} A_{k+1}^{j}{\scriptstyle (j,k+1)}  ((U_l^{u} z_{k+1}^{j} + U_{k+1}^{j} z_l^{u} - U_{k+1}^{j} U_{l}^{u})) \nonumber\\
 \end{numcases}
\end{scriptsize}

To sum up, we relax Constraints~(\ref{ReluRelaxed}) as follows:
  \begin{numcases}{}
   z_{k+1}^{j} (z_{k+1}^{j} - \displaystyle{\sum_{u=1:~u \text{ unstable}}^{n_k}} (W_{k+1,u}^{j} z_k^{u})) - b_{k+1}^j)  \nonumber\\\quad\quad\;-B{\scriptstyle (j,k+1)}  z_{k+1}^{j}\leq UP_r & $\forall r=1,2 $ \ref{RelaxReluUP}\\
    z_{k+1}^{j} (z_{k+1}^{j} - \displaystyle{\sum_{u=1:~u \text{ unstable}}^{n_k}} (W_{k+1,u}^{j} z_k^{u})) - b_{k+1}^j)  \nonumber\\\quad\quad\; -B{\scriptstyle (j,k+1)}  z_{k+1}^{j}\geq LO_r & $\forall r=1,2$ \ref{RelaxReluLO}
  \end{numcases}

Note that a wider set of constraints could be explored by varying combinations of bounds from~(\ref{stab:bound1}) and~(\ref{stab:bound2}). This ablation multiply by 4 the number of ReLU  constraints dedicated to unstable neurons. However, it eliminates several categories of constraints related to stable active neurons, including the ReLU  constraint~(\ref{ReLURef1_SDP_k}), the bounding constraints~(\ref{boundRef_SDP_k}), the triangular constraints~(\ref{tri_SDP}), the McCormick constraints~(\ref{mc4-obj}), and the RLT constraints~(\ref{RLT_p_SDP}). Notably, this leads to a quadratic reduction with respect to $n_k^{a}$ :~$4~n_k^{a}~|\setJKbar|$~constraints removed from~\ref{mc4-obj}, and up to $n_k^{a} (n_k^{u}+n_k^{a})$ constraints removed from~(\ref{RLT_p_SDP}) on layer $k$. In sufficiently big neural networks, this quadratic reduction counterbalances the linear increase in constraints with respect to $n_k^{u}$.

\subsubsection{Representation on a small example}
We write here an example of this method on a small neural network, with $K = 4$ layers, whose hidden layers have 3 neurons each : $n_1 = n_2 = n_3 = 3$. The input and output are bidimensional $n_0 = n_4 = 2$. We specify here the weights of the DNN : 

\begin{equation*}
    \begin{array}{llll|llll}
        W_1 &= \begin{pmatrix} 1 & 2 & 3 \\ 0 & 2 & -1\end{pmatrix} & b_1 &= \begin{pmatrix}1 &-2 & 1 \end{pmatrix} & W_2 &=  \begin{pmatrix} 1 & 0 & 3 \\ -2 & 2 &  -1 \\ 4& -1 & 0 \end{pmatrix} & b_2 &= \begin{pmatrix} -2 & 1 & -1 \end{pmatrix}    \\[5ex]
        W_3 &= \begin{pmatrix} 0 & -2 & 3 \\ 1 & 2 & -1 \\ -1 & 1 & -2 \end{pmatrix} &b_3 &= \begin{pmatrix}3 & -3 & 5 \end{pmatrix} & W_4 &=  \begin{pmatrix} 0 & -2 \\ 4 & 0 \\ 0 & 1 \end{pmatrix} & b_4 &= \begin{pmatrix} 2 & -4  \end{pmatrix}   \\
    \end{array}
\end{equation*}

Suppose that neurons $1$ and $2$ of layer $1$, $0$ of layer $2$, $1$ of layer $3$ are stable active. Neurons $1$ of layer $2$ and $0$ of layer $3$ are inactive. Formally, if we represent the given ReLU neural network with stable active neurons in green and stable inactive neurons in red, we have : 

\begin{center}
\begin{tikzpicture}[scale=2,      transform shape]

        \node (I1) at (0,0.8) [circle, draw=black, fill=black!50] {};
        \node (I2) at (0,0.3) [circle, draw=black, fill=black!50] {};
      
        \node[align=left,anchor=west, font = \fontsize{6}{8}\selectfont] at (-0.2,0.8) {\color{white} 0};
        \node[align=left,anchor=west, font = \fontsize{6}{8}\selectfont] at (-0.2,0.3) {\color{white} 1};
        
        \node (H11) at (1,1) [circle, draw=black, fill=black!50] {};
        \node (H12) at (1,0.5) [circle, draw=black, fill=green!50] {};
        \node (H13) at (1,0) [circle, draw=black, fill=green!50] {};

         \node[align=left,anchor=west, font = \fontsize{6}{8}\selectfont] at (0.8,1) {\color{white} 0};
        \node[align=left,anchor=west, font = \fontsize{6}{8}\selectfont] at (0.8,0.5) {\color{white} 1};
         \node[align=left,anchor=west, font = \fontsize{6}{8}\selectfont] at (0.8,0) {\color{white} 2};

        \node (H21) at (2.3,1) [circle, draw=black, fill=green!50] {};
        \node (H22) at (2.3,0.5) [circle, draw=black, fill=red!50] {};
        \node (H23) at (2.3,0) [circle, draw=black, fill=black!50] {};

            \node[align=left,anchor=west, font = \fontsize{6}{8}\selectfont] at (2.1,1) {\color{white} 0};
         \node[align=left,anchor=west, font = \fontsize{6}{8}\selectfont] at (2.1,0.5) {\color{white} 1};
         \node[align=left,anchor=west, font = \fontsize{6}{8}\selectfont] at (2.1,0) {\color{white} 2};

        \node (H31) at (3.5,1) [circle, draw=black, fill=red!50] {};
        \node (H32) at (3.5,0.5) [circle, draw=black, fill=green!50] {};
        \node (H33) at (3.5,0) [circle, draw=black, fill=black!50] {};

        \node[align=left,anchor=west, font = \fontsize{6}{8}\selectfont] at (3.3,1) {\color{white} 0};
        \node[align=left,anchor=west, font = \fontsize{6}{8}\selectfont] at (3.3,0.5) {\color{white} 1};
         \node[align=left,anchor=west, font = \fontsize{6}{8}\selectfont] at (3.3,0) {\color{white} 2};
         
        \node (O1) at (4.5,0.8) [circle, draw=black, fill=black!50, font = \fontsize{6}{8}\selectfont] {};
        \node (O2) at (4.5,0.3) [circle, draw=black, fill=black!50, font = \fontsize{6}{8}\selectfont] {};

        \node[align=left,anchor=west, font = \fontsize{6}{8}\selectfont] at (4.3,0.3) {\color{white} 1};
        \node[align=left,anchor=west, font = \fontsize{6}{8}\selectfont] at (4.3,0.8) {\color{white} -3};
        
        \foreach \i in {1,2}
            \foreach \j in {11,12,13}
                \draw (I\i) -- (H\j);

        \foreach \i in {11,12,13}
            \foreach \j in {21,22,23}
            \draw (H\i) -- (H\j);
        \foreach \i in {21,22,23}
            \foreach \j in {31,32,33}
                \draw (H\i) -- (H\j);
        \foreach \i in {31,32,33}
            \foreach \j in {1,2}
                \draw (H\i) -- (O\j);
    \end{tikzpicture}
    \end{center}

Let us compute the linear decomposition of each stable active neuron.
\begin{itemize}
    \item $z_1^{1} = W_1^{1} z_0 + b_1^{1} = 2 \times z_0^{0} + 2 \times z_0^{1} - 2$, which gives $\lambda_0^{0}{\scriptstyle (1,1)} = 2$, $\lambda_0^{1}{\scriptstyle (1,1)} = 2$, $\gamma{\scriptstyle (1,1)} = - 2$.
    \item $z_1^{2} = W_1^{1} z_0 + b_1^{1} = 3 \times z_0^{0} - 1 \times z_0^{1} + 1$, which gives $\lambda_0^{0}{\scriptstyle (2,1)} = 3$, $\lambda_0^{1}{\scriptstyle (2,1)} = -1$, $\gamma{\scriptstyle (2,1)} = 1$.
    \item We have $\lambda_1^{0}{\scriptstyle (0,2)} = 1$, $\lambda_0^{0}{\scriptstyle (0,2)} = 8$, $\lambda_0^{1}{\scriptstyle (0,2)} = -8$, $\gamma{\scriptstyle (0,2)} = 6$ thanks to the following calculus:
    \begin{align*}
        z_2^{0} &= W_2^{0} z_1 + b_2^{0} \\
        &= 1 \times z_1^{0} - 2 \times z_1^{1} + 4 \times z_1^{2} - 2 \\
        &= z_1^{0} -2 \; (2 z_0^{0} + 2 z_0^{1} - 2) + 4 \; (3 z_0^{0} - z_0^{1} + 1) - 2 \\
        &= z_1^{0} + 8 z_0^{0} - 8 z_0^{1} + 6 \\
    \end{align*}
    \item We have $\lambda_2^{2}{\scriptstyle (1,3)} = 1$, $\lambda_1^{0}{\scriptstyle (1,3)} = -2$, $\lambda_0^{0}{\scriptstyle (1,3)} = -16$,  $\lambda_0^{1}{\scriptstyle (1,3)} = 16$ $\gamma{\scriptstyle (1,3)} = -15$ thanks to the following calculus:
     \begin{align*}
        z_3^{1} &= W_3^{1} z_2 + b_3^{1} \\
        &= -2 \times z_2^{0} + 2 \times z_2^{1} + 1 \times z_2^{2} - 3 \\
        &= -2 \big( z_1^{0} + 8 z_0^{0} - 8 z_0^{1} + 6 \big) + 2 \times 0 + z_2^{2} - 3   \text{  Neuron 1 of layer 2 is inactive} \\
        &= z_2^{2} - 2 z_1^{0} -16 z_0^{0} +16 z_0^{1}  - 15    \\
    \end{align*}
\end{itemize}

We represent now the quadratic ReLU constraint on the unstable neuron $2$ of layer $3$. Formally, the exact constraint is : 

$$z_3^{2} \times \big(z_3^{2} - (3 z_2^{0} - z_2^{1} - 2 z_2^{2} + 5) \big) = 0$$

where neuron $1$ of layer $2$ is inactive and can be removed. Separating the stable active neuron $0$ and unstable neuron $2$ of layer $2$, the constraint gives: 
\begin{equation}
    \begin{array}{ll}
        z_3^{2} \times \big(z_3^{2}  + 2 z_2^{2} - 5 \big) &=  3 z_2^{0} z_3^{2}  \\
        &= 3 \big(  z_1^{0} + 8 z_0^{0} - 8 z_0^{1} + 6 \big) z_3^{2} \\
        &=  3 z_1^{0} z_3^{2} + 24 z_0^{0} z_3^{2} - 24 z_0^{1} z_3^{2}+ 18 z_3^{2} \\
    \end{array}
    \nonumber
\end{equation}

We have $B{\scriptstyle (2,3)} = 18$, $A_0^{0}{\scriptstyle (2,3)} = 24, A_0^{1}{\scriptstyle (2,3)} = -24, A_1^{0}{\scriptstyle (2,3)} = 3$.
We consider now the following bounding given by the McCormick constraints : 

\begin{equation*}
\begin{cases}
\begin{array}{ll}
    3 z_1^{0} z_3^{2} &\in [0,~~~3~U_1^{0}~z_3^{2}] \\
    3 z_1^{0} z_3^{2} &\in [3~U_3^{2}~z_1^{0},~~~3~U_3^{2}~z_1^{0}+~3~U_1^{0}~z_3^{2}-3~U_{3}^{2}~U_1^{0}] \\[5ex]

    24 z_0^{0} z_3^{2} &\in [24~\tilde{L}_0^{0}~z_3^{2},~~~24~U_0^{0}~z_3^{2}] \\
    24 z_0^{0} z_3^{2} &\in [24~U_3^{2}z_0^{0}+24~\tilde{L}_0^{0}z_3^{2}-24~\tilde{L}_0^{0}U_3^{2},24~U_3^{2}z_0^{0}+24~U_0^{0}z_3^{2}-24~U_{3}^{2}U_0^{0}] \\[5ex]

    -24 z_0^{1} z_3^{2} &\in [-24~U_0^{1}~z_3^{2},~~~-24~\tilde{L}_0^{1}~z_3^{2}] \\
    - 24 z_0^{1} z_3^{2} &\in [-24~U_3^{2}z_0^{1}-24~U_0^{1}z_3^{2}+24~U_{3}^{2}U_0^{1},-24~U_3^{2}z_0^{1}-24~\tilde{L}_0^{1}z_3^{2}+24~\tilde{L}_0^{1}U_3^{2}]
\end{array}
\end{cases}
\end{equation*}

By summing up the first of the two boundings, we obtain the following inequalities : 
\begin{equation}
    \begin{array}{lll}
      C_2~z_3^{2} &\leq z_3^{2} \times \big(z_3^{2}  + 2 z_2^{2} - 5 \big) &\leq  C_1~z_3^{2}\\
    \end{array} \nonumber
\end{equation}

with $C1 = 3~U_1^{0}~+24~U_0^{0}-24~\tilde{L}_0^{1}$ and $C_2 = 24~\tilde{L}_0^{0}-24~U_0^{1}$. \\
By summing up the second of the two boundings, we obtain the following inequalities : 
\begin{equation*}
    \begin{array}{l}
      C_{4,3}^{2} z_3^{2} + C_{4,1}^{0} z_1^{0} + C_{4,0}^{1} z_0^{1} + C_{4,0}^{0} z_0^{0} + B z_3^{2} + C_4  \\\quad\leq z_3^{2} \times \big(z_3^{2}  + 2 z_2^{2} - 5 \big) \quad \leq C_{3,3}^{2} z_3^{2} + C_{3,1}^{0} z_1^{0} + C_{3,0}^{1} z_0^{1} + C_{3,0}^{0} z_0^{0} + B z_3^{2} +  C_3 \\
    \end{array}
\end{equation*}

with $C_{3,3}^{2} = 3 U_1^{0} + 24 U_0^{0} - 24 \tilde{L}_0^{1}$,~$C_{3,1}^{0} = 3 U_3^{2}$,~$C_{3,0}^{0} =  24 U_3^{2}$,~$C_{3,0}^{1} = -24 U_3^{2}$, $C_{3} = U_3^{2} (-3 U_1^{0} - 24 U_0^{0} + 24 \tilde{L}_0^{1})$, and $C_{4,3}^{2} = 24 \tilde{L}_0^{0} - 24 U_0^{1}$,~$C_{4,1}^{0} = 3 U_3^{2}$,~$C_{4,0}^{0} =  24 U_3^{2}$,~$C_{4,0}^{1} = -24 U_3^{2}$, $C_{4} = 24~U_3^{2} (U_0^{1} - \tilde{L}_0^{0})$,

\subsection{Triangular constraint}
\label{supp_seq:pp_triang}
We say that a neuron $j$ of layer $k$ of lower bound $L_k^{j}$ and upper bound $U_k^{j}$ is \emph{stable active} if $L_k^{j} \geq 0$ and \emph{stable inactive} if $U_k^{j} \leq 0$. Otherwise, a neuron is \emph{unstable} (i.e. if $L_k^{j} < 0 < U_k^{j}$). \\

In order to tighten the upper bound on the ReLU  activation function, a well-known constraint is the triangular constraint, which provides a convex embedding of the ReLU  output. Depending on the activation status of the considered neuron $j$ of layer $k$, we decompose the linear upper bound in the following set $\setT$ of triangular inequalities:

\begin{numcases}{z^j_{k} \in \setT \Leftrightarrow}
    z^j_{k} \leq 0 & if $j$ is \emph{inactive} \nonumber \\
     z^j_{k} \leq W_k^{j} z_{k-1} + b_k^{j} & if$j$ is \emph{ active} \nonumber \\
      z^j_{k} \leq \frac{U_k^j}{U_k^j - L_k^j} (W_k^{j} z_{k-1} + b_k^{j}) \nonumber \\ 
      \qquad + \frac{U_k^j}{U_k^j - L_k^j}(b_k^j - L_k^j)  & if $j$ is \emph{unstable} \nonumber 
\end{numcases}

Combined with constraints $z_k^{j} \geq 0$, $z_k^{j} \leq \hat{z}_k^{j}$, where $\hat{z}_k^{j}$ denotes preactivation vector of neuron $j$ of layer $k$, this constraint yields the exact output $z_k^{j} = 0$ when neuron $j$ is inactive, and $z_k^{j} = W_k^{j} z_{k-1} + b_k^{j}$ when it is active.

\begin{figure}[h]
    \centering
    \includegraphics[width=0.55\linewidth]{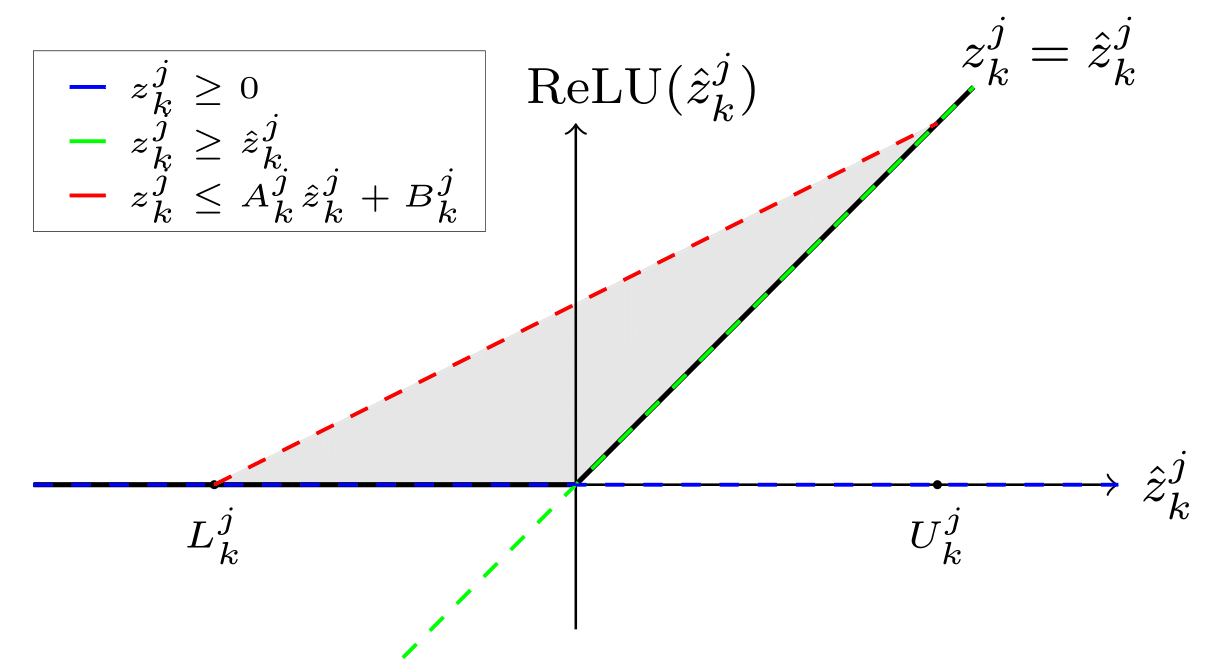}
    \caption{Triangular constraint on neuron $j$ of layer $k$, where $\hat{z}_k^{j} = W_k^{j} z_{k-1} + b_k^{j}$ is the preactivation vector}
    \label{fig:triangular_constraint}
\end{figure}

When neuron $j$ is unstable, the upper-bound is plotted in red in Figure~ \ref{fig:triangular_constraint}.
    The equation of the red line clearly depends on the lower and upper bounds $L^j_k$ and $U^j_k$ of the preactivation vector $W_k^{j} z_{k-1} + b_k^{j}$. 
    \begin{small}
   \end{small}

The triangular constraint has been shown to be limited, as it represents only the convex hull of the output of a single ReLU  neuron~\citep{salman_provably_2019}. Recent works have proposed convex relaxations to capture the joint behavior of multiple ReLUs. Such ideas could be explored to compute more efficient cuts in SDP models.

\subsection{Coherence constraint~\cref{stab_SDP}}
In all our SDP models, we use the constraint between two consecutive matrices $P_k[(1 \; z_{k+1})(1 \;z_{k+1})^{\top}] = P_{k+1}[(1 \; z_{k+1})(1 \; z_{k+1})^{\top}]$~(\ref{stab_SDP}) relaxed as in~\citep{batten_efficient_2021} and in~\citep{lan_tight_2022}. For a layer $k$, including all coherence constraints in the model would introduce $\frac{n_k (n_k +3)}{2}$ constraints. In the relaxation of constraint~(\ref{stab_SDP}), only the $n_k$ linear constraints remain: $P_k[z_{k+1}]=P_{k+1}[z_{k+1}]$, preventing the number of constraints from exploding. 

\subsection{RLT constraint}
We present here the RLT cuts~(\ref{RLT_p_SDP}) selected in~\citep{lan_tight_2022}.
\label{supp_seq:app_rlt}
These cuts contribute to tightening the relaxation, and are given below:
\begin{numcases}{\text{Tightening cuts}}
  \quad  z_{k}z_{k+1} \leq  L_kz_{k+1} + U_{k+1}z_{k} - L_k U_{k+1}  \label{RLT_1} \\
   \quad z_{k}z_{k+1} \leq  U_k z_{k+1} + L_{k+1}z_{k} - U_k L_{k+1} \label{RLT_2} \\
   \quad z_{k}z_{k+1} \geq  L_k z_{k+1} + L_{k+1}z_{k} - L_kL_{k+1} \label{RLT_3} \\
   \quad P_k[z_{k+1}z_{k+1}] \leq  U_{k+1} z_{k+1} \label{triangular_RLT_1}\\
   \quad P_k[z_{k+1}z_{k+1}^T] \leq  A_k P_k[z_k z_{k+1}] + B_k P_k[z_{k+1}] \label{triangular_RLT_2}
\end{numcases}

The number of these constraints is large. For a given layer $k$, constraints~(\ref{RLT_1})--(\ref{RLT_3}) scale quadratically with $n_{k+1} \times n_k$, while constraints~(\ref{triangular_RLT_1}) and~(\ref{triangular_RLT_2}) scale with $n_{k+1} \times n_{k+1}$. Including all these constraints in the SDP model would significantly increase the computing time. A heuristic is therefore needed to select only a subset of these cuts. As~(\ref{triangular_RLT_1}) and~(\ref{triangular_RLT_2}) capture \textit{intra} layer dependencies, a heuristic selecting a subset of them is difficult to design. In contrast, since~(\ref{RLT_1})--(\ref{RLT_3})  represents \textit{inter} layer dependencies, a heuristic based on the linear layer weights linking them is possible. \\

Only a subset of constraints~(\ref{RLT_1})--(\ref{RLT_3}) is finally selected, based  on a given percentage $p$. Specifically, for each neuron $j$ on layer $k+1$ we select $\lfloor p \, n_{k} \rfloor$ cuts. The heuristic sorts the absolute value of the weights $|W_{k+1}^{j}|$, and selects the neurons corresponding to the top $\lfloor p \, n_{k} \rfloor$ entries in the sorted vector. This selection is based on the full size of layer $k$ regardless of whether neurons have been pruned. More precisely, in the context of a full ablation, we do not select $\lfloor p \, n_{k}^{u} \rfloor$  RLT cuts but $\lfloor p \, n_{k} \rfloor$ ones. 

\subsection{Tightening cuts for SDP$_U$}

For simplicity and clarity, we used the full logits $z_K^{j}$ in constraints~(\ref{eq:cut1})~(\ref{eq:cut2})~(\ref{mc4-obj}). Note that these logits are not variables of our model. To obtain the full constraints in our model, we need to substitute each logit by its linear expression with respect to the penultimate layer variables : $z_K^{j} = W_K^{j} z_{K-1} + b_{K}^{j}$. 

Furthermore, for the sake of clarity in constraints~(\ref{mc1})~(\ref{mc4-obj})~(\ref{eq:cut1})~(\ref{eq:cut2}~(\ref{cut:betaibetaj}), we omit explicit matrix indexation. To recover the full constraint, one must, for example, replace variables such as $\beta_i \beta_j$ by $P_{K-2}[\beta_i \beta_j]$. Combining this with the logit expression, $z_K^{j}$ should be replaced by $W_K^{j} P_{K-2}[z_{K-1}] + b_{K}^{j}$, and $\beta_j z_K^{j}$ by $W_K^{j} P_{K-2}[\beta_j z_{K-1}] + b_{K}^{j} P_{K-2}[\beta_j]$.

\newpage
\section{Experimentations}

We ran all our experiments on the MOSEK solver, with a number of threads of 4 and a time limit for a single optimization of $7200$s.

\subsection{Implementation details}
\label{supp_seq:exp}
All networks have been trained with a batch size of 128, the Adam optimizer, and a learning rate of 0.001. For reproducibility, we show the details of the adversarial training in table \ref{tab:benchmark_implementation_details}.  All PGD attacks were used with number of steps $= 40$, a random start. We denote by \texttt{6x100-5}, \texttt{6x100-20}, \texttt{6x100-50}, \texttt{6x100-67} networks used in experiment 3.

\begin{table*}[h]
\centering
\renewcommand{\arraystretch}{1.3}
\begin{tabular}{|l|l|c|l|c|c|}
    \hline
    \textbf{Network} & \textbf{Architecture} & \multicolumn{2}{c|}{\textbf{Adversarial Training}} & \multicolumn{2}{c|}{\textbf{Accuracy (\%)}} \\
    \cline{3-6}
     & & \textbf{Epochs} & \textbf{Adversarial attack} & \textbf{Clean} & \textbf{PGD} \\
    \hline
    \texttt{9x100} & 784-9x100-10 & 200 & PGD ($\epsilon = 0.3$, $\alpha = 0.01$) & 95.2 & 95 \\
    \texttt{9x200} & 784-9x200-10 & 200 & PGD ($\epsilon = 0.3$, $\alpha = 0.01$) & 96.9 & 100 \\
    \hline
    \texttt{FCNNA} & 3072-2x20-100 & 200 & PGD ($\epsilon = \frac{8}{255}$, $\alpha = 0.005$) & 13.4 & 13.6 \\
    \hline
    \texttt{6x100-5}  & 784-6x100-5  & 100 & PGD ($\epsilon = 0.3$, $\alpha = 0.01$) & 97.6 & 100 \\
    \texttt{6x100-20} & 784-6x100-20 & 100 & PGD ($\epsilon = 0.3$, $\alpha = 0.01$) & 86.0 & 80  \\
    \texttt{6x100-50} & 784-6x100-50 & 100 & PGD ($\epsilon = 0.3$, $\alpha = 0.01$) & 76.7 & 80  \\
    \texttt{6x100-67} & 784-6x100-67 & 100 & PGD ($\epsilon = 0.3$, $\alpha = 0.01$) & 75.2 & 82  \\[1ex]
    \hline
\end{tabular}
\caption{Networks used in our three experiments. Clean accuracy refers to the accuracy on the validation set. PGD accuracy denotes the empirical percentage of robustness against the PGD attack (overestimation of the actual robustness) on the subdataset  evaluated on a subdataset whose size varies across experiments : 100 samples from \texttt{MNIST} for \texttt{9x100} and \texttt{9x200}; 1000 samples from \texttt{CIFAR100} for \texttt{FCNNA}; 5 samples from \texttt{MNIST} for \texttt{6x100-5}; 20 samples from \texttt{MNIST} and \texttt{FashionMNIST} for \texttt{6x100-20}; 50 and 67 samples from \texttt{MNIST}, \texttt{FashionMNIST}, \texttt{EMNIST} and \texttt{KMNIST} for \texttt{6x100-50} and  \texttt{6x100-67} respectively.}
\label{tab:benchmark_implementation_details}
\end{table*}

\newpage
\section{Additional state of the art}

\subsection{Complete verification}
\label{supp_seq:complete_verification}

Ideal verification is complete, ensuring that all answers are reliable. However, due to the complexity of the problem, fully achieving such verification is often constrained. 
Works using Satisfiability Modulo Theory~\citep{ehlers_formal_2017} like ReLUplex~\citep{katz_reluplex_2017} or Marabou~\citep{gurfinkel_marabou_2024} have been developed, but are not currently scalable. Nevertheless, they are very precise and give formal proof of robustness or useful counterexamples when working on a sufficiently small network. 

Some works have introduced Mixed Integer Programming formulations (MIP), see ~\citep{fischetti_deep_2018,dsouza_maximum_2017},
but the direct resolution of these models without relaxations is also not scalable. Some works propose an efficient selection of the variables that should be relaxed~\citep{liao_solution-aware_2025}. Most efficient complete verification methods rely on Branch \& Bound~\citep{bunel_branch_2020,ferrari_complete_2022,jaeckle_neural_2021,lu_neural_2019}, whose relaxation of the certification problem is fast heuristics like CROWN. They have been improved by smart branching: splitting on the activation or not activation of a set of neurons has been exponentially faster than splitting on the input ball. 
Some methods have improved ReLU  splitting to better choose neurons for branching decision~\citep{henriksen_deepsplit_2021, zhou_scalable_2024}.  When reaching a certain depth of the tree, a relatively fast MIP is solved (with few binary variables as the activation of most neurons is fixed), to prune a branch without exploring all its content. It also helps to avoid impossible activation patterns in practice, which may not be seen by heuristic methods like bound propagation, guaranteeing the soundness of the algorithm. 
The resolution of these MILP has been further improved with cutting planes in method GCP-CROWN~\citep{zhang_general_2022}.  Other works are Branch\&Bound frameworks based on SDP relaxations ~\citep{lan_semidefinite_2023,chiu_sdp-crown_2025}.  Note that these approaches primarily improve certification rates rather than scalability.

\subsection{Adversarial training}
\label{supp_seq:app_createDNN}

In this section, we present the adversarial training used in order to create robust networks. Madry introduced adversarial training \citep{madry_towards_2019} by adding a maximization problem into the common training minimization problem:

$$\min\limits_{\forall (x,y) \in \mathcal{X}} \max\limits_{z \in \mathcal{B}_{\epsilon}(x)} \mathcal{L}(z,y) $$

where the inner maximization represents the computation of the worst adversarial attack. It is approached by heuristics, most of them are based on gradient descents, Projected Gradient Descent (PGD)~\citep{madry_towards_2019}, or its variants (FGSM  \citep{goodfellow_explaining_2015}, IGS~\citep{kurakin_adversarial_2017}, JSMA~\citep{papernot_limitations_2015} for the norm $\|.\|_0$, DeepFool~\citep{moosavi-dezfooli_deepfool_2016} for $\|.\|_2$) and have been improved since ~\citep{carlini_towards_2017}. 
Specifically, we train our model using the most classic adversarial attack: PGD, and compute untargeted adversarial attacks with it. \\

\bibliographystyle{splncs04}
\bibliography{../Certification_Problem_bibtex}
\nocite{Cor76, ehlers_formal_2017, katz_reluplex_2017, gurfinkel_marabou_2024, fischetti_deep_2018, dsouza_maximum_2017, liao_solution-aware_2025}
\nocite{bunel_branch_2020,ferrari_complete_2022,jaeckle_neural_2021,lu_neural_2019, henriksen_deepsplit_2021, zhou_scalable_2024, zhang_general_2022, lan_semidefinite_2023,chiu_sdp-crown_2025}
\nocite{madry_towards_2019, goodfellow_explaining_2015, kurakin_adversarial_2017, papernot_limitations_2015, moosavi-dezfooli_deepfool_2016, carlini_towards_2017, salman_provably_2019}

\nocite{goodfellow_explaining_2015, katz_reluplex_2017, tjeng_evaluating_2019, fischetti_deep_2018, dsouza_maximum_2017, wong_provable_2018}
\nocite{raghunathan_semidenite_2018, wang_beta-crown_2021, Ans09, vandenberghe_chordal_2015, batten_efficient_2021, ehlers_formal_2017, lan_tight_2022, SheAda90}
\nocite{Cor76, jung2020continual, tjeng_evaluating_2019, serra_lossless_2020, botoeva_efficient_2020}
\nocite{mosek, wang_beta-crown_2021, raghunathan_semidenite_2018, batten_efficient_2021, lan_tight_2022, li_sok_2023, singh_abstract_2019}
\nocite{fischetti_deep_2018, tjeng_evaluating_2019, ehlers_formal_2017, katz_reluplex_2017, wong_provable_2018, Gowal_ijcai2019, weng_towards_2018, huang2017safety, raghunathan_semidenite_2018, zhang_tightness_2020, dathathri_enabling_2018, chiu_tight_2023, lan_iteratively_2023, batten_efficient_2021, lan_tight_2022, han_accelerating_2024}